# New probabilistic methods for physics

Guy Cirier

LSTA, University Pierre et Marie Curie Sorbonne, France
Email: guy.cirier@gmail.com

**Résumé**
On présente une méthode probabiliste pour étudier les itérations. Cette méthode est utilisée pour des EDO et est appliquée pour analyser les comportements asymptotiques en mécanique ou en physique.

**Abstract**
We present how a probabilistic model can describe the asymptotic behavior of the iterations, with applications for ODE and approach of some problems in mechanics in $\mathbb{R}^d$.

## Introduction

Since a long time, probabilistic methods are used in physics. They concern generally a great quantity of objects having the same behavior (particles, atoms, stars…). Then statistical methods can be applied as the number of objects is very great. These results are well known.
So, Boltzmann has developed the concept of ergodism for study the kinetic theory of the gaz. Roughly speaking, in mathematics, the ergodic hypothesis consists to write that the time average equals spatial mean. Then, a probability of presence can be defined in the space. Many theoretical results exist, but the hypothesis is often open to debates.

But the notion of probability of presence can be defined intrinsically, without ergodism. We consider an invariant measure $P$ under a measurable function $f$:
*$P$ is invariant under $f$ if, for all borelian set $B$, $P$ verifies the Perron-Frobenius's equation PF* (See Lasota & Mc Kay[6]):

$$P_f(B) = P\circ f^{-1}(B) = P(B).$$

This equation doesn't suppose that the borelians are bounded. However, with this choice, many demonstrations are simplified and it assumes easily the convergence of the studied series. Effectively, if $f$ is iterated indefinitely in a bounded or compact set, the probability of presence gets a sense with the more or less important density of points in some places. If such a $P$ can be defined explicitly, we get important information about the objects verifying the $f$ iteration. But this functional equation is not easy to study. However, it's the subject of the present paper.

So, $f$ is supposed be a function applying a bounded set $C \subset \mathbb{R}^d$ in itself. This condition is verified for all phenomena of the infinitely small or for many phenomena on the earth. Let $\Delta \subset C$ the support of this invariant measure associated. Naturally, iterations of $f$ means $f^{(k)} = f \circ f \circ \ldots \circ f$ for every $k$. Let the measure $P$ invariant: $P_f(B) = P\circ f^{-(k)}(B) = P_{f^{(k)}}(B)$.

Let $\Delta_k$ be the domain of $P_{f^{(k)}}$: $\Delta_k \subset \Delta_{k-1} \subset \Delta \subset C$, because $B \subset f^{-1} \circ f(B)$. If $f$ is injective, we have: $\Delta_k = \Delta_{k-1} = \Delta$. Implicitly, we suppose that the iteration starts from a point in $C$. But, a function $f$ may have many invariant domains.

**Main results concerning physics**:
This method brings new explications for many bounded phenomena:
- Under general conditions, this measure exists, in particular if $f$ is $C^\infty$ and applies a bounded set $C \subset \mathbb{R}^d$ in itself. It implies self-similarity in some cases.
- We consider a differential ordinary equation $\mathrm{d}a/\mathrm{d}t = F(a)$ as an iteration $f$: $f(a) = a + \delta F(a)$ where $\delta = t/n$ is the path. if $f$ is $C^\infty$ and applies a bounded $C \subset \mathbb{R}^d$ in itself, the final behavior is asymptotically almost periodic or the dimension of the phenomena is flattening.
- But, there are some levels hard to cross continuously due to the Fredholm's alternatives. At that time, breakings of the behavior can appear. The Hamilton's equation is reviewed. The n-bodies problem has new perspective.

These results will be clarified in the following paper presented in two parts with chapters:
**Part I: Invariant measure of Perron-Frobenius**
**Part II: Differential ordinary equations as iterations**

# Part I: Invariant measure of Perron-Frobenius

## A- Laplace-Fourier's transform of the invariant measure *P*

### 1- Laplace-Fourier's transform of the invariant measure *P*

**Hypothesis**
*$C \subset \mathbb{R}^d$ is bounded. $f$ applies $C$ in $C$ and $f$ at least $C^\infty$. Then this measure $P$ exists.*
Here, we search an analytic approach of the $f$-invariant measure $P$ with the Fourier-Laplace's transform. We use the known property of this measure: for all positive $P$- measurable function $g$, we have the formula (see Lasota & Mc Kay[6]):

$$\int g \circ f(x)\, dP(x) = \int g(x) dP(x).$$

For $g(x) = e^{yx}$, we write the Fourier-Laplace's transform $\emptyset(y) = \mathcal{L}(e^{yX}) = E\left(e^{yf(X)}\right)$ . As $x \in C \subset \mathbb{R}^d$ is bounded, $|x| < D = diam(C)$, we develop $\emptyset(y)$ with series $\emptyset(y) = \Sigma_n b_n y^n$. It is convergent because $|b_n| < E(\frac{|x|^n}{n!}) < D^n/n!$ . Let $\emptyset_f(y) = \mathcal{L}\left(e^{yf(X)}\right)$. The second series is also convergent because $|f(x)| < D$. If the measure $P$ is invariant:

$$\emptyset(y) = \emptyset_f(y).$$

Here, $y \in \mathbb{R}^d$ ou $\mathbb{C}^d$. If $y = it$, $\emptyset(y)$ is the characteristic function of the measure $P$. In this case, we recall that equality or convergence of the characteristic functions implies equality or convergence in law. This will be often implied in this paper.

***Translation of the distribution with a small fixed vector*** $a \in \mathbb{R}^d$, $X \longmapsto X + a$.

We translate $X$ : $\quad \emptyset(y, a) = E\left(e^{y(X+a)}\right)$.

And by $f$-transformation [1]): $\quad \emptyset(y, a) \longmapsto \emptyset_f(y, a) = E\left(e^{yf(X+a)}\right)$.

As the measure is invariant: $\quad \emptyset(y, a) = \emptyset_f(y, a)$.

And: $\quad \theta(y, a) = \emptyset(y, a) - \emptyset_f(y, a) = 0$.

***Proposition:***
*$\theta(y, a) \equiv 0$ is an identity and $\partial^p \partial^q \theta(y, a)/\partial y^p \partial a^q \equiv 0$ for all $p$ and $q$ . For $\forall a \in C$ and $\forall y$:*

$$\theta(y, a) = \Sigma_n b_n \partial^n (e^{ya} - e^{yf(a)})/\partial a^n \equiv 0.$$

■ For $a = 0$, $\theta(y, 0) \equiv 0$. For $\forall a \neq 0$, $e^{-ya} \theta_f(y, a) = \emptyset(y, 0) - E\left(e^{yf(X+a) - ya}\right) = 0$. So, $\theta(y, a) \equiv 0$ is an identity.

If the random vector $X\epsilon C \subset \mathbb{R}^d$ has a measure $P$ with density $p(x)$, the translated random vector $X+a$ has the same density for every small translation $a \in C \subset \mathbb{R}^d$ of the random vector $X$. Using the convergent series $\emptyset(y) = \Sigma_n b_n y^n$, we have the translated density[2]:

$$p(x-a) = \mathcal{L}^{-1}\left(e^{ta}\,\emptyset(t)\right) = \mathcal{L}^{-1}(\,\Sigma_n b_n t^n e^{ta}).$$

So, we can write $p(x-a)$ as a distribution in the sense of Schwartz with the Dirac's $\delta$:

$$p(x-a) = \Sigma_n b_n \partial^n \delta(x-a)/\partial a^n.$$

As: $$E\left(e^{y(X+a)}\right) = \Sigma_n b_n y^n e^{ya} = \Sigma_n b_n \partial^n e^{ya}/\partial a^n.$$

And : $$E\left(e^{yf(X+a)}\right) = \int e^{yf(x+a)}\,dP(x) = \int e^{yf(x)}\,p(x-a)\,dx\,,$$

$$E\left(e^{yf(X+a)}\right) = \Sigma_n b_n \partial^n (\int e^{yf(x)}\,\delta(x-a)dx)/\partial a^n,$$

$$E\left(e^{yf(X+a)}\right) = \Sigma_n b_n \partial^n e^{yf(a)}/\partial a^n.$$

By difference, we get: $$\theta(y,a) = \Sigma_n b_n \partial^n (e^{ya} - e^{yf(a)}\,)/\partial a^n \equiv 0\,.\blacksquare$$

**Notations**

*We call* $\theta(y,a) \equiv 0$ *resolving equation* $R_a$ *of PF and* $e^n(y,a) = \partial^n(e^{ya} - e^{yf(a)}\,)/\partial a^n$ *gap of order* $n$*:*

$$\theta(y,a) = \Sigma_n b_n e^n(y,a) \equiv 0,$$

*If* $a=0$*:* $$\theta(y) = \theta(y,0) = \Sigma_n b_n e^n(y,0) \equiv 0.$$

**Remarks**

- We observe that $\partial^n(e^{yf(a)}\,)/\partial a^n = H_n(y,a)\,e^{yf(a)}$where $H_n(y,a)$ is a Bell-polynomial in $y$ with degree $n$ even $f$ is $C^\infty$or analytic. We can note the gap:

$$e^n(y,a) = \partial^n(e^{ya} - e^{yf(a)}\,)/\partial a^n = y^n e^{ya} - H_n(y,a)\,e^{yf(a)}$$

And, for $a=0$: $e^n(y) = e^n(y,0) = y^n - H_n(y)$ is a polynomial in $y$ with degree $n$.

- Writing $\emptyset(y) = \Sigma_n b_n y^n$, we obtain $\emptyset_f(y)$ by putting $H_n(y)$ instead of $y^n$ in the series.

- If $a=0$, then: $\emptyset(0) = 1,\ \emptyset_f(0) = 1 = e^{yf(0)}, \forall y$ and $f(0)=0$. The point of reference is a fixed point. If we take $f^{(k)}$ instead of $f$, the point of reference, yet denoted 0, is a point of a cycle. These points, and the eigen values $\lambda$ of the linear part of $f$ at these points, are most significant for the study of the convergence of the process. They are well known for the linearization of $f$. That means: find a function $\varphi$ such $\varphi \circ f = \lambda_1 \varphi$. The problem is linked to the eigen values $\lambda$ in the resonance's case ${\lambda_1}^n$=1, not studied here.

- If $y=0$, then: $\emptyset_f(0) = \emptyset(0) = 1$ and $b_0 = 1$. But, the other $b_n$ are unknown.

**2- Iteration induces derivation on** $\theta(y,0) = 0$

**Proposition**

*If* $\lambda_\ell \neq 1$, $\ell = 1,\ldots, d$, *the iteration* $f_\ell(a)$ *acts as a derivation on* $\theta(y,0) = 0$ *in the sense*:

$$a_\ell \longmapsto f_\ell(a) \Longrightarrow\ \theta(y,0) \longmapsto \partial\theta(y,a)/\partial a_\ell|_{a=0}\ ,$$

*denoted:* $$a \longmapsto f(a) \Longrightarrow \theta(y,0) \longmapsto \partial\theta(y,a)/\partial a|_{a=0}.$$

*By induction, all the coordinates of* $n \in N^d$ *in* $e^n(y,0)$ *are:* $n_1 = \cdots = n_\ell \ldots = n_d$ *without resonances:*

$$n = \cdots = n_\ell \ldots = n_d.$$

■ For all derivable function $g(a)$ in $\mathbb{R}$, such as $g(0)=0$, we write $g(a) \sim a\partial g(a)/\partial a$ near 0. Let a function $f(a)$ such as $f(0)=0$ and $\partial f(a)/\partial a|_{a=0} = \lambda \neq 1$ at 0. Then, we have:

$$gf(a) \sim a\lambda\,\partial g(a)/\partial a|_{a=0}.$$

We apply this result to the impact of the transformation $a_\ell \longmapsto f_\ell(a)$ on $D = \theta\left(y, f_\ell(a)\right) - \theta(y,a_\ell) = 0$ , as all the other coordinates $\overline{a_\ell} \neq a_\ell$ remain fixed:

As: $$\theta(y,a_\ell) = 0\,,$$

when $a_\ell \to 0$: $$f_\ell(a) - a_\ell \sim a_\ell(\lambda_\ell - 1)$$

$$D \sim a_\ell(\lambda_\ell - 1)(\theta(y,0)/\partial a_\ell).$$

As an iteration of $f$ is the iteration of each coordinate:

$$D \sim \Pi_{\ell=1}^{\ell=d} a_\ell(\lambda_\ell - 1)(\partial\theta(y,0)/\partial a).\blacksquare$$

We remark that the iteration $f_\ell(a)$ acts also as a derivation on the gap $e^n(y) = 0$.

**3- The asymptotic lattice distribution of probability**

**Proposition**

*The general solution of the linear equation* $\theta(y) = 0$ *is* $\emptyset(y) = 1 + b\varphi(y)$ *with an arbitrary constant b. It means that* $\emptyset(y) = 1$, *for all* $\varphi(y) = 0$. *We have a lattice distribution of probability for* $\varphi(y) = 0$.

■ $\emptyset(y)$ can be written $\emptyset(y) = 1 + b\varphi(y)$ with $\varphi(0) = 0$. The linearity of $\theta_f(y) \equiv 0$ implies: $b\varphi(y) - b\varphi_f(y) = b\left(\varphi(y) - \varphi_f(y)\right) \equiv 0, for\ \forall b$ . So: $\emptyset(y) = 1 + b\varphi(y), \forall b$. ■

**4-To summarize this paragraph,** without resonances, we can write the Laplace- Fourier transform $\emptyset(y)$ of the invariant measure: $\emptyset(y)=1 + b\varphi(y) = \Sigma_m b_m y^m$, $\emptyset_f(y) = 1 + b\varphi_f(y) = \Sigma_m b_m H_n(y)$ and: $\theta(y) = \emptyset(y) - \emptyset_f(y) = \Sigma_{1\le m} b_m(y^m - H_m(y))$. Then:

- As $b$ is arbitrary, $\varphi(y) = 0$, $\varphi_f(y) = 0$ and $\theta(y) = \emptyset(y) - \emptyset_f(y) \equiv 0$.
- As each iteration induces a derivation on $\theta(y) = 0$, we can limit some studies to $n \in N^d$ , with $n = n_1 = \ldots = n_\ell \ldots = n_d$, for instance, we can study $\varphi_n(y) = \Sigma_{m\le n} b_m y^m$ with a common $n$ , $\varphi_{nf}(y) = \Sigma_{m\le n} b_m H_m(y)$ and: $\theta_n(y) = \Sigma_{m\le n} b_m(y^m - H_m(y))$.
- As series are convergent, we have: $|\varphi_n(y) - \varphi(y)| < \epsilon$, then: $|\varphi_n(y)| < \epsilon$; simultaneously: $\left|\varphi_{nf}(y)\right| < \epsilon$, $|\theta_n(y)| < \epsilon$ uniformly for all bounded $y \in \mathbb{R}^d$ with non-null coefficients.

## B- Solution of the Perron-Frobenius's equation

$n \in N^d$, but we have seen that we can take a maximal common index: $n = n_1 = \ldots = n_\ell \ldots = n_d$.
Let $h_{mk}$ the known coefficients of the Bell's polynomial of $f$: $H_m(y) = \Sigma_{0<k\le m} h_{mk}\, y^k$.

***Lemma***

*We write the polynomial:* $\theta_n(y) = \varphi_n(y) - \varphi_{nf}(y) = \Sigma_{m\le n} b_m(y^m - H_m(y))$ *according to the powers of* $y^m$*:*

$$\theta_n(y) = (1-\lambda^n) b_n\, y^n + \Sigma_{m<n}(b_m - \Sigma_{m\le j\le n} h_{jm}\, b_j - b_n h_{nm})y^m$$

■ The Bell's polynomial is: $H_m(y) = \Sigma_{0<k\le m} h_{mk}\, y^k$. We note:

$$\Theta_n(y) = (1-\lambda^n) b_n\, y^n - b_n \Sigma_{m<n} h_{nm}\, y^m + \Sigma_{m<n} b_m(y^m - \Sigma_{0<k\le m} h_{mk}\, y^k).$$

Then, we rearrange the expression with the powers of $y^m$. ■

As $\theta_f(y) \equiv 0$ , we must search an estimator $\theta *_n (y)$ of the polynomial $|\theta_n(y)| < \epsilon$ , verifying also $|\varphi_n(y)| < \epsilon$, $\left|\varphi_{nf}(y)\right| < \epsilon$. A non-null $\theta *_n$ must verify the following conditions: $|\theta *_n (y)| < \epsilon$, $|\varphi *_n (y)| < \epsilon$ and $\left|\varphi *_{nf} (y)\right| < \epsilon$. For fixed $n$, such estimator can be obtained with a stricter condition: $\theta *_n (y) = 0$ ; but, we have now a polynomial which must be null. Let: $Y_n = \{y | e^n(y) = 0\}$ and $Y *_n = \{y | \theta *_n (y) = 0 = 0\}$. We will see that $y \in Y_n$ induces $y \in Y *_n$.

For fixed $n$, rewrite $\varphi *_n (y) = \Sigma_{m\le n} b *_m\, y^m$, $\varphi *_n (y) = \Sigma_{m\le n} b *_m\, H_m(y)$ and $\theta *_n (y) = \Sigma_{m\le n} b *_m\, e^m(y) = 0$. For $y = 0$: $b *_0 = 1$ and $e^0(y) = 0$. The other $b *_m$ are unknown.

We have to find a sequence of $b *_m$, $m < n$ such as $\theta *_n (y) = 0$, we get the solution.

If $\theta *_n (y) = 0$, without resonances, $(1-\lambda_1{}^n) b_n\, y^n$ cannot be null except if $e^n(y) = 0$. We put the arbitrary constant $b = b *_n$.

***Theorem***

*Without resonances, we find asymptotically a unique convergent solution of* $\theta *_n (y) = 0$, *up to an arbitrary constant b:*

$$\emptyset_n(y) = 1 - be^n(y).$$

*We obtain a lattice distribution defined by the zeros of* $e^n(y)$.

*In the repellent case where* ${\lambda_1}^n \gg 1$, *we have:*

$$\emptyset_n(y) \sim 1 - bH_n(y).$$

*Then, the distribution of the real zeros of the polynomials* $H_n(y)$ *gives the distribution of the Perron-Frobenius's measure when* $n \to \infty$.

*We obtain a lattice distribution defined by the zeros of* $H_n(y)$.

■ - By definition: $\theta *_n (y) = \Sigma_{m\leq n} b *_m e^m(y) = 0$ for an arbitrarily fixed $b *_n = b \neq 0$.

$$\theta *_n (y) = be^n(y) + \Sigma_{m<n} b *_m e^m(y) = 0,$$

$$\Sigma_{m<n} b *_m e^m(y) = -be^n(y),$$

With the lemma, the polynomial $\theta *_n (y)$ is null without resonances:

$$(1\text{-}{\lambda_1}^n) b\, y^n + \Sigma_{m<n} \left(b *_m - \Sigma_{m\leq j\leq n} h_{jm}\, b *_j - bh_{nm}\right) y^m = 0.$$

For $\forall\ y^m$ such as $m < n$, if $b *_m - \Sigma_{m\leq j\leq n} h_{jm}\, b *_j - bh_{nm} = 0$, we obtain a finite triangular system of linear equations: it is easy to solve step by step and all the unknown coefficients $b *_m$ defined in a unique way in function of $b$ and the known coefficients $h_{mk}$ of $H_m(y)$ with $m \leq n \in N^d$ . This is thru for $\forall y$, but now $(1\text{-}\lambda{\lambda_1}^n) b\, y^n \neq 0$ except if we chose: $y \in Y_n = \{y | e^n(y) = 0\}$. For these $y$ , we have simultaneously:

$\theta *_n (y) = be^n(y) + \Sigma_{m<n} b *_m e^m(y) = 0$ and $be^n(y) = 0$.

So, for $\forall m < n$, we have: $b\Sigma_{0<k<n} h_{nk}\, y^k + \Sigma_{m<n} b *_m e^m(y) = 0$.

Conversely, if $e^n(y) = 0$,it easy to see that $b *_m$are solution. This non-null solution is unique for $\forall b *_n = b \neq 0$ fixed.

- We must prove that the $b *_m$ converge to the $b_m$. As we can write $|\varphi_n(y) = \alpha(y)| < \epsilon$, with $\alpha(y) = \Sigma_{m<n} \alpha_m y^m$, for $y$ belonging to $Y_n$, then we replace $e^n(y)$ by $e^n(y) + \alpha(y)$ and the equations become:

$$b_m - \Sigma_{m\leq j\leq n} h_{jm}\, b_j - b_n h_{nm} + \alpha_m = 0.$$

As the $\alpha_m \to 0$, by continuity of the linear equations with fixed coefficients $h_{nm}$: $b *_m \to b_m$.

So, we have constructed the polynomials $\emptyset *_n (y) - 1$ and $\emptyset_f *_n (y) - 1$ under the condition no-depending of $n$, and we can write $\emptyset *_n (y) = 1 + b\varphi *_n (y) = 1 + b\varphi *_{nf} (y)$ where $b$ is arbitrary. We simplify the notation with writing: $\emptyset_n(y) = 1 - be^n(y)$.

Different cases can happen according to ${\lambda_1}^n \gg 1$ or ${\lambda_1}^n \ll 1$. For instance, if all the coordinates of $|\lambda_1|$are less than 1, the process converges to the fixed point. If some of them are less than 1, but others are greater than 1, we have a hyperbolic situation without resonances.

When $y^n {\lambda_1}^n \gg y^n$ , we can write for large $n$ :

$$\emptyset_n(y) \sim 1 - bH_n(y). ■$$

And now we have to study the zeros of $H_n(y)$.

## C- Asymptotic density of the real zeros of $H_n(y)$

### 1- Study of the asymptotic density of the real zeros of $H_n(y)$

We study now the real zeros of $H_n(y)$ in order to get their asymptotic distribution. $f$ is real, $C^\infty$, without resonances, and applies $C \subset \mathbb{R}^d$ bounded in itself. As each iteration induces a derivation on $H_n(y) = 0$, we take in general a common index $n = n_1 = n_\ell \ldots = n_d$.

We recall that the polynomial:

$$H_n(y) = e^{-yf(a)}\partial^n e^{yf(a)}/\partial a^n|_{a=0} = \partial^n e^{yf(a)}/\partial a^n|_{a=0} ,$$

can be represented by the Cauchy's integral:

$$H_{n-1}(y) = c\oint_\Gamma \frac{e^{yf(a)}}{a^n} da = c\oint_\Gamma e^{yf(a)-n\ln a} da,$$

where $\Gamma$ is a closed polydisk around $a = 0$ of $f$, $a \in \mathbb{C}^d$, $c$ is finite non-null real, without importance in the context. We take $n = n_1 = n_\ell = n_d$.

The integrand is: $$n\gamma(a) = yf(a) - n\ln(a).$$

If we note $y = ns$ ($s_\ell = y_\ell/n$): $$\gamma(a) = sf(a) - \ln(a).$$

$\gamma(a)$ is called here the Plancherel-Rotach's function [8] because they are the firsts to use this method to approximate the Hermite's polynomials. The method consists to search the critical point $\alpha$ of $\gamma(a)$. Under the numerous conditions of the general position**,** the critical point $\alpha$ maximizing $e^{n\gamma(a)}$ gives an approximation of the integral. This critical point $\alpha$ is defined by the equation:

$$\partial\gamma(\alpha)/\partial a = s\partial f(\alpha)/\partial a - 1/\alpha = 0.$$

The nullify of the approximation gives us the real zeros of $H_n(y)$. Some coordinates of $\alpha$ can be real, the others are complex. $H_\infty(s) = 0$ only if $\alpha$ have complex coordinates.

**Main conditions to use this method**

The critical point $\alpha$ must be isolated from the other critical points, at a finite distance and far from the edge $\partial C$. A sufficient condition to get this maximum is that the hessian matrix of $\gamma(a)$, which is hermitian, is definite negative at $\alpha$. that means $x[\partial^2 yf(a)/\partial a^2]\, x < 0$. But, as the Morse's theory says that the hessian matrix has the form $\Sigma' x'_{\ell'}{}^2 - \Sigma\, x_\ell{}^2$, the positive coordinates are not suitable. Moreover, when the hessian is not definite negative, we have many other problems with linearities, Stokes's phenomena or catastrophe theory. In other hand, the method is valid with few modifications for integrals with boundaries [5]. But, what happens if other critical points appear ?

## 2- The asymptotic density for $f$ real

### Theorem

*Under the conditions of the general position, the critical point* $\alpha$ *of the PR-function must be complex to annul* $H_n(y)$. *Among the multiple critical points, first, we must choose* $\alpha$ *such as the real part of the Cauchy's integral is maximum.*

*- Then, only for all complex coordinates (yet denoted* $\alpha$ *or* $\Im\mathrm{m}()$*) of* $\alpha$*, we have:*

$$\partial\gamma(\alpha)/\partial a = s\partial f(\alpha)/\partial a - 1/\alpha = 0.$$

$$n\Im\mathrm{m}(\gamma) = \Im\mathrm{m}(nsf(\alpha) - n\ln\alpha) = k\pi.$$

*When* $n \to \infty$, $k/n = \chi_n \to \chi$, *with* $\chi$ *uniform on* $[0,1]$ *for each r complex coordinates.*

*The asymptotic distribution is:*

$$\Im m(sf(\alpha) - \ln\alpha) = \chi\pi. \quad \text{With } \kappa \in [0,1]^r,$$

*and:* $$\Im m(\gamma_\ell(\alpha))) - \chi_\ell\pi = 0.$$

*The* $\kappa_\ell$ *are identically independent uniform on* $[0,1]$.

*- In the unidimensional case, the distribution of the zeros is:*

$$q(s)\,ds = \Im m\big(f(\alpha)\big)\,ds/\pi \text{ with } s\partial f(\alpha)/\partial a - 1/\alpha = 0.$$

■ For all real coordinates of $\alpha$, the exponential defining $H_{n-1}(y)$ cannot be cancelled. Only a complex coordinate can annul $H_{n-1}(y)$. Suppose there are *r* complex coordinates. So, as $f$ has real coefficients, if $\partial\gamma(\alpha)/\partial a = 0$; then, $\partial\gamma(\overline{\alpha})/\partial a = 0$ where $\overline{\alpha}$ is the conjugate of $\alpha$. As in Plancherel- Rotach, the contribution to estimate $H_{n-1}(y)$ is:

$$H_{n-1}(y) = c\oint_\Gamma \left(\exp\big(n\gamma(a)\big) - \exp(n\gamma(\overline{a}))\right) da.$$

- As in the steepest descent, $\Im\mathrm{m}(n\gamma)$ is constant:

$$H_{n-1}(y) = 2c\sin\,(\Im\mathrm{m}(n\gamma)\,\mathcal{R}e(\oint_{\Gamma} e^{n\gamma(a)}\mathrm{d}a).$$

For the maximum of the real part, only the imaginary part $\Im\mathrm{m}(\gamma)$ of $\gamma(a)$ can nullify $H_{n-1}(y)$. That means:

$$\sin\,(\Im\mathrm{m}(n\gamma)) = 0 \text{ and } \Im\mathrm{m}(n\gamma) = k\pi,$$

$$\Im\mathrm{m}(sf(\alpha)) - \theta = \pi k/n = \pi\chi_n,$$

Writing: $\gamma_\ell(\alpha)) = s_\ell f_\ell(\alpha) - \ln(\alpha_\ell)$, : $\gamma(\alpha) = \Sigma_\ell \gamma_\ell(\alpha))$ and with: $\chi = \Sigma_\ell \chi_\ell$,

$$\Sigma_\ell (n\Im\mathrm{m}(\gamma_\ell(\alpha) - \chi_\ell\pi) = 0.$$

But, at time $n + 1_\ell$, where all the $n$ except $n_\ell = n$ which becomes $n_\ell + 1_\ell$, the PR function becomes $\Sigma_\ell (n\Im\mathrm{m}(\gamma_\ell(\alpha) - \chi_\ell\pi) + \gamma_\ell(\alpha) - \chi_\ell\pi = 0.$

Because $H_{n-1}(y) = 0$ becomes $H_{n-1+1_\ell}(y) = 0.$

then: $\Im\mathrm{m}(\gamma_\ell(\alpha))) - \chi_\ell\pi/n = 0.$

When $n \rightarrow \infty$, $k_\ell/n \rightarrow \chi_\ell \in [0,1]$,

with $\kappa_\ell$ uniform on [0,1] for each complex coordinate and we have $H_\infty(s)$=0.

- In the case of unidimensional function, the repartition of the zeros verifies for $\chi = k/n$:

$$\Im\mathrm{m}\big(\gamma(\alpha)\big) = \chi\pi \quad \text{with} \quad \partial\gamma/\partial a = 0.$$

So: $q\,(s)\,\mathrm{d}s$ =Prob {1 zero $\in (s, s + \mathrm{d}s)$} $=\mathrm{d}\chi$,

$$q\,(s)\,\mathrm{d}s = \Im\mathrm{m}(\mathrm{d}\gamma/\mathrm{d}s)\,\mathrm{d}s/\pi = \Im\mathrm{m}\big(f(\alpha)\big)\,ds/\pi\,,$$

because: $\mathrm{d}\gamma/\mathrm{d}s = \partial\gamma/\partial s + \partial\gamma/\partial a.\partial a/\partial s = f(\alpha)$. ■

**Remarks**

- As some coordinates of $\alpha$ can be real, then stochastic manifolds may appear.
- These distributions of the PF-equation are linked to each fixed (or cycle) point $f(0) = 0$. Then, the solutions are local. They preexist independently of all convergence, as the fixed points, cycles or invariant varieties. The convergence to one of the asymptotic situations depends on the basin of the initial point.
- All these distributions can be masked in various situations. The principle of the maximum of the real part $\Re e(\gamma)$ of $\gamma(a)$ provides a method to define the fuzzy frontiers of the different basins of attraction.
- In the case where $H_{n-1}(ns) = 0$, an iteration induces a derivation as we have seen for $\theta_f(y, 0) = 0$ . The relation $k/n < \frac{(k+1)}{n+1} < (k + 1)/n$ can explain that the zeros of $H_n$ are interleaved when the reciprocal image $s$ of $\Im\mathrm{m}\big(sf(\alpha)\big) - \theta = k\pi/n$ are monotonic.
- If $f$ is complex, we can identify $\mathbb{C}$ as $\mathbb{R}^2$ and the calculus are quite the same but more difficult.

For instance, if $f$ has real coefficients, with $s$ real, $\alpha$ and $\kappa$ fix: $\gamma(a) = sf(a) - \ln\,(a\,)$. The distribution of zeros is defined by:

$$\Im\mathrm{m}\big(\gamma(\alpha)\big) = k\pi/n\ ;\ \partial\gamma(\alpha)/\partial\alpha = s\partial f(\alpha)/\partial\alpha - 1/\alpha = 0.$$

Let $\xi = e + id$ a complex number, we want to study the iterations of $\xi f$. The equations of the distribution are now defined for the PR function at the critical point $\alpha'$:

$$\gamma'(\alpha') = s'\xi f(\alpha') - \ln(\alpha')$$

$$\Im\mathrm{m}(\gamma') = k'\pi/n \text{ and } \partial\gamma(\alpha')/\partial\alpha' = s'\xi\partial f(\alpha')/\partial\alpha' - 1/\alpha' = 0.$$

The steepest descent induces: $H_{n-1}(y) = 2c\sin\,(\Im\mathrm{m}(n\gamma')\,\mathcal{R}e(\oint_{\Gamma} e^{n\gamma'(a)}\mathrm{d}a.$

- Most important is the fact that iterations are asymptotically commutative as derivations.
- As everybody knows, the steepest descent's method is difficult to use, but it shows a very large variety of behaviors.

**3- Relation between the distribution of the zeros of $H_n(y)$ and the invariant density**

We calculate the distribution $P$ with the asymptotic distribution of zeros of $H_n(y)$ in a particular case. We simplify in supposing that we have $n^d$ points $y_n \in \mathbb{R}^d$ verifying $H_n(y) = 0$. The Laplace's transform of the $y_n$ is:

$$\emptyset_n(y) = (1/n^d)\Sigma_n e^{-yy_n}).$$

At each point $y_n$, correspond a point $x$ in $C$ with non-null density. But, we must normalize the coordinates of the zeros $y$ to get a distribution at finite distance: $nx_{n\ell} = y_{n\ell} = ns_\ell$ . At the volume $dy$ correspond a normalized volume $dx$: $dy = n^d dx$. So, we get $n^d q(x)dx$ points.
When $n \to \infty$, we have:

$$\emptyset_n(ns) = (1/n^d)\Sigma_n e^{-sx_n} n^d q(x_n)dx_n ,$$

$$\emptyset_n(ns) \to \int e^{-sx}\, q(x)dx = \emptyset(\mathrm{s}).$$

$\emptyset_n(\mathrm{s})$ is the Laplace's transform of the density of zeros of $H_n(y)$.
$\emptyset(\mathrm{s})$ is asymptotically the Laplace's transform of the invariant distribution $q(x)$ the density of zeros solution of $H_n(y) = 0$. Let $p(x)$ be the density of the measure invariant with regard to the Lebesgue's measure of $\mathbb{R}^d$ ($dP(x) = p(x)dx$). What is the relation of $p(x)$ with the asymptotic density of zeros $q(x)$? Hint:

**Proposition**
*Suppose that C is rectangular with null density* $p(x)$ *on the boundaries* $\partial C$. *Then* $p(x)$ *is:*

$$p(x) = (-x\partial q\,(x)/\partial x).$$

■ For all solution of $\theta\ (y, a) = 0$, the measure $P$ must verify the identity $\theta\ (y, a) \equiv 0$, for $\forall y$ and
small $\forall a$. If $\varphi(y)$ is solution of $\theta\ (y, a) = 0$, then, the Laplace's transform must verify:

$$\varphi(y) = \partial^2\emptyset(\mathrm{y}, \mathrm{a}) \,/\, \partial a \partial y|_{a=0} = \partial(y\emptyset(y))/\partial y.$$

The inverse transformation gets the result. In fact, there are few phenomena where we have $n^d$ points $y \in \mathbb{R}^d$ verifying $H_n(y) = 0$ , except the case unidimensional because all the zeros of $H_n(y)$ are not real. ■

**Remark: the method of the kernel for the invariant measure**
This method seems complicated, but more productive than others more elegant but difficult to use as the following method of the kernel:
Let $\emptyset(y)$ the Fourier-Laplace's transform of the invariant measure $P$ corresponding to $f$:

$$\emptyset(y) = \emptyset_f(y),$$

$$\emptyset_f(y) = E\left(e^{yf(X)}\right) = \mathcal{L}\left(e^{yf(X)}\right) = \int e^{yf(x)}\, dP(x).$$

But: $$p(x) = \mathcal{L}^{-1}\left(\,\emptyset(t)\right) = \frac{1}{(2\pi i)^d}\int_{c-i\infty}^{c+i\infty} e^{-tx}\emptyset(t)dt.$$

So: $$\emptyset_f(y) = \int e^{yf(x)}\, dP(x) = \frac{1}{(2\pi i)^d}\ \int e^{yf(x)} \int_{c-i\infty}^{c+i\infty} e^{-tx}\emptyset(t)dtdx.$$

As $\emptyset(y) = \emptyset_f(y)$ : $$\emptyset(y) = \frac{1}{(2\pi i)^d}\int\int_{c-i\infty}^{c+i\infty} e^{yf(x)-tx}\,\emptyset(t)dtdx.$$

We recognize the Legendre's transformation in the integrand:

$$C(y,t) = \frac{1}{(2\pi i)^d}\int_{c-i\infty}^{c+i\infty} e^{yf(x)-tx}dx.$$

We obtain the functional: $\emptyset(y) = \int C(y,t)\, \emptyset(t)dt,$
where $C(y,t)$ is the kernel for the eigen value 1. The steepest descent can be used to approximate $C(y,t)$, but it doesn't seem very productive.

**4- Analysis of the domination in the case of many fixed points**
When we have many fixed points, there is a problem to know the domain of attraction of each point. Let $\alpha$ and $\beta$ two points of $Fix(f)$ in a bounded set $C \subset \mathbb{R}^d$. $f$ applies $C \subset \mathbb{R}^d$ bounded in itself. Let $a \in C$: $a = \alpha + u = \beta + v$,
Changing the origin: $f_\alpha(u) = f(\alpha + u)$- $\alpha$ and $f_\beta(v) = f(\beta + v)$- $\beta$,

with $f(\alpha)$- $\alpha$= $f(\beta)$- $\beta$=0. Let: $\gamma_\alpha\ (u) = sf_\alpha(u) - \ln\ (u)$ and: $\gamma_\beta\ (u) = sf_\beta(v) - \ln\ (v)$.
The critical point is defined by: $s\partial f(\alpha + u)/\partial u - 1/u = 0$ et $s\partial f(\beta + v)/\partial v - 1/v = 0$,
as: $a = \alpha + u = \beta + v$: $\quad v = u$,
The steepest descent method that we must take the greatest real part of the integrand at the critical point. By difference, the sign of the real part of $s(f_\alpha(u) - f_\beta(u))$ at the critical point get the criterium. But this frontier may be stochastic.

**5- Examples**

- *Let the logistic map* [9]: $f(a) = \lambda a - a^2/2$; and $\gamma(a)$=$s\ (\lambda a - a^2/2) - \ln(\,a)$;

$$\partial\gamma/\partial a = s(\lambda a - a^2) - 1 = 0\,.$$

If the discriminant: $\quad \Delta = (s\lambda)^2 - 4s < 0$ and $s > 0$,
we put $\lambda\sqrt{s} = 2cos\vartheta$ , we have: $2\cos\vartheta\, a\sqrt{s} - sa^2 - 1 = 0$ with roots: $a\sqrt{s} = e^{\pm i\vartheta}$ ;
and: $\quad \Im\mathrm{m}(f(a)) = \Im\mathrm{m}(\lambda a - a^2/2) = \sin 2\vartheta/s$.

$$q\ (s)\ \mathrm{d}s = \Im\mathrm{m}(f(a))ds/\pi = (1 - \cos 2\vartheta)\mathrm{d}\vartheta/\pi,$$

So: $\quad q\ (s) = (\lambda/2\pi)\sqrt{1/s - \lambda^2/4}\,.$

If we put: $\quad t = \cos\vartheta = \frac{\lambda\sqrt{s}}{2}$ ,

Then $t$ follows: $\quad W(t)\ \mathrm{d}t = (2/\pi)\ \sqrt{1 - t^2}\mathrm{d}t.$

We recover directly a well-known result: Let $H_n\ (t, a)$= $\partial^n(e^{t(\lambda a - a^2/2)})/\partial a^n$ generating the Hermite polynomials $H_n\ (t)$. The law of the zeros of $H_n(t)$ is known as the semi-circular Wigner's law: $W(t)\mathrm{d}t = (2/\pi)\sqrt{1 - t^2}\mathrm{d}t$.
- Then, the density of the logistic corresponding to $q\ (s)$ is $p(a)$:

$$p(a) = -a\mathrm{d}q/\mathrm{d}a = -a(2/\pi)\mathrm{d}\Big(\frac{\lambda}{4\sqrt{a}}\sqrt{1 - \frac{a\lambda^2}{4}}\Big)/\mathrm{d}a,$$

$$p(a) = \lambda/(2\pi\sqrt{4a - a^2\lambda^2}).$$

The density of the logistic map follows a Beta (1/2,1/2) low in a more general situation than in the Ulam-Von Neumann's case [9].
- *m-Hermitian case*: $\quad f(a) = \lambda a - a^m/m.$
The Plancherel-Rotach's function is: $\gamma(a)$=$s(\lambda a - a^m/m) - \ln(a)$.
With the critical point $a$ defined by the trinomial equation:

$$\mathrm{d}\gamma(a)/\mathrm{d}a = s(\lambda a - a^m) - 1 = 0,$$

studied by H. Fell. We note that the discriminant $\Delta < 0$ may establish the existence of a complex root of this equation. We recall: $\quad \Delta = m^m/s^{m-1} - (m-1)^{m-1}\lambda^m < 0.$
- *Quadratic function*: We take now a quadratic function $f$ in $\mathbb{R}^d$ with $f(0) = 0$. We write the PR function $\gamma(a)$ for every fixed point 0 of $f(a)$:

$$f(a) = \lambda a + Qa^2/2,$$

the hessian of $sQ$ is symmetric. For all $s$ such as $sQ$ is non -degenerate, it exists an orthogonal transformation $T$: $a = Tu$, with: $T'sQT = D$, the diagonal matrix of eigenvalues of $sQ$ and:

$$\ln\ (a) = \Sigma_{\ell=1}^{\ell=d}\ln(a_\ell) = \ln\Pi_{\ell=1}^{\ell=d}(a_\ell) = \ln\ (\mathrm{Vol}(a))$$
$$= \ln\ (\mathrm{Vol(u)}) = \ln\ (\Pi_{\ell=1}^{\ell=d}u_\ell) = \ln\ (u) = \Sigma_{\ell=1}^{\ell=d}\ln(u_\ell),$$

because the volume $\mathrm{Vol}(a) = \mathrm{Vol(u)}$ is invariant under an orthogonal transformation.
Then, the P.R. function $\gamma(a)$ becomes:

$$\gamma(u) = sf(Tu) - \ln Tu = s\lambda Tu + Du^2/2 - \ln u.$$

We note $D = D^+{}_\ell$ if $\ell = 1, .., p$ and $D = D^-{}_\ell$ if $\ell = p + 1, .., d$. Then:

$$\gamma(u) = \Sigma_{\ell=1}^{\ell=d}\Lambda_\ell u_\ell + \Sigma_{\ell=1}^{\ell=p}D^+{}_\ell {u_\ell}^2/2 - \Sigma_{\ell=p+1}^{\ell=d}D^-{}_\ell {u_\ell}^2/2 - \Sigma_{\ell=1}^{\ell=d}\ln u_\ell,$$
$$= \Sigma_{\ell=1}^{\ell=p}(\Lambda_\ell u_\ell + D^+{}_\ell {u_\ell}^2/2 - \ln u_\ell) + \Sigma_{\ell=p+1}^{\ell=d}(\Lambda_\ell u_\ell + D^-{}_\ell/2{u_\ell}^2 - \ln u_\ell),$$

where: $\quad \Lambda u = s\lambda Tu.$

If we note: $\gamma_+(u_\ell) = \Lambda_\ell u_\ell + D^+{}_\ell {u_\ell}^2/2\text{-ln}u_\ell$ ,
and $\gamma_-(u_\ell) = \Lambda_\ell u_\ell + D^-{}_\ell {u_\ell}^2/2\text{-ln}u_\ell$,
So: $\gamma(u) = \Sigma_{\ell=1}^{\ell=p}\gamma_+(u_\ell) + \Sigma_{\ell=p+1}^{\ell=d}\gamma_-(u_\ell)$.
Applying the logistic calculus to each $\gamma_+(u_\ell)$ and $\gamma_-(u_\ell)$, we obtain $p$ conditions $\Lambda_\ell u_\ell = 0$ half the time and *d-p* random independent variables following a Beta (1/2,1/2) low. But, we may have other fixed points: $a(1-\lambda) = Qa^2$.
These results can be extended to a $C^\infty$ function $f$ with the Morse-Palais Lemma as if the hessian is definite. We have studied $f$, we examine now $h = g \circ f$ and $f^{(p)}$.

## D- Self similarity

### 1- Analysis of $h = g \circ f$ and applications to $f^{(p)}$

Let $f$ and $g$ two functions $C^\infty$ which apply a bounded set $C \subset \mathbb{R}^d$ in itself. We search invariant measures of $h = g \circ f$ from the study of the invariant measures of $f$.
- We know $(g \circ f)^{-1} = f^{-1} \circ g^{-1}$. If the measure $P$ is invariant under $f$, $P \circ f^{-1}(B) = P(B)$,
So: $P((g \circ f)^{-1}(B)) = P(f^{-1} \circ g^{-1}(B)) = P(g^{-1}(B))$.
- We have to search zeros of $H_{n-1,h}(y)$ for $h = g \circ f$ in relation with those of $H_{n-1,f}(y)$ for $f$. We suppose that the fixed point 0 of $f$ is fixed point of $g$: $g(0) = 0$ .
For $f$, let α, $s_f$ solution of $H_{n-1,h}(y) = 0$ completely defined for fixed $\chi$ by:
$\Im m\ (\gamma_f(\alpha\ ) = \chi\ \pi$ and $\partial\gamma_f/\partial a = 0$.
We have to search solution ($s_{g\circ f}$ et $\alpha_{g\circ f}$) for $h = g \circ f$.

**Lemma 2**

*Under the previous hypothesis, if we choose* $s_{g\circ f}$ *verifying* $\partial \frac{g(a)}{\partial a}|_{a=f\,(\alpha)} = s_f$, *the critical complex point* $\alpha$ *of* $\gamma_f$: $\alpha_f = \alpha$ *remains the critical point* $\alpha_{g\circ f} = \alpha$ *for the iteration* $g \circ f$. *Then, the distribution of* $h = g \circ f$ *provided by the reciprocal image of the distribution of* $f$ :
$\Im m\ (\gamma_{g\circ f}\ (\alpha\ ) = \chi_{g\circ f}\,\pi$.
■ Let: $\gamma_f\ (a) = s_f f(a) - \ln a$ and $\gamma_{g\circ f}\ (a) = s_{g\circ f}\ g \circ f(a) - \ln a$.
The critical point $\alpha$ for $f$ is defined by:
$\partial\gamma_f/\partial a = 0$,
and: $s_f \partial f(\alpha)/\partial a = 1/\ \alpha$.
The critical point $\alpha_{g\circ f}$ for $g \circ f$ is defined by:
$\partial\gamma_{g\circ f}/\partial a = 0$,
That means: $s_{g\circ f}\partial\ g \circ f(\alpha_{g\circ f})/\partial a = 1/\ \alpha_{g\circ f}$.
But: $\partial g \circ f\ (\alpha_{g\circ f})/\partial a = \partial g(f)/\partial f\,.\,\partial f(\alpha_{g\circ f})/\partial a$.
If we take: $s_{g\circ f}\,\partial\ g(f)/\partial\, f = s_f$,
Then: $s_f \partial f(\alpha_{g\circ f})/\partial a = 1/\ \alpha_{g\circ f}$ .
The critical point $\alpha_{g\circ f}$ verifies the critical equation of $f$, so $\alpha_{g\circ f} = \alpha$. ■

**Remarks**

Even if $s_f$ is real, $s_{g\circ f}$ is not real because $\partial g(f)/\partial f$ depends on $\alpha$. It is not a big problem because the steepest descent is valid for a complex function. But the complex part of $s_{g\circ f}$ can also be nullify. In order to have $\partial\ g(f)/\partial\ f$ independent from $a$, $\partial^2 g(f)\partial a^2\ \ \partial f(a)/\partial a$ =0. The hessian of $g$ must be degenerated at the critical point $\alpha$. However, $g$ may induce other solutions and linearities.

**Corollary**

*g= $f^{(p)}$, then the critical point $\alpha$ of $\gamma_f$ is critical point of $\gamma_{f^{(p)}}$ , $\forall\ p$.*

**2- Analysis of $f^{(p)}$, $\forall\ p$**

Some distributions can be linked to cycles of order $p$. By translation the origin at a point $0 \in \mathbb{R}^d$of the cycle of $f^{(p)}$: $f^{(p)}(0)$ =0, the iteration is written:

$$h(a)= f^{(p)}(a) - f^{(p)}(0) )= f^{(p)}(a) )= f^{(p-1)} \circ f\ (a).$$

But, before all, we have to study the behavior of these distributions near a fixed realization of the asymptotic random vector $\chi$ uniform on $[0,1]^{\ r}$, This means that $\partial\gamma(\alpha)/\partial a = 0$ has complex solutions.

**Lemma 3**

*For$\forall\ \mu \in \mathbb{Z}$, the asymptotic distribution of zeros of $H_n(\mu y)$ is the same than this of $H_n(y)$and asymptotically, if $\mu \in \mathbb{R}$, we have the same result.*

■ Let $\kappa$ a realization $\chi$ and a fixed $\alpha$ such as $H_\infty(s) = 0$: $\partial\gamma(\alpha)/\partial a = s\partial f(\alpha)/\partial a - 1/\alpha$ =0.
If $\mu$ is rational, $\mu = p/q$ with $p, q \in \mathbb{N}$, then: $n\partial\gamma(\alpha)/\partial a = y\partial f(\alpha)/\partial a - n/\alpha$ =0.
Now $n\partial\gamma_\mu(\alpha)/\partial a = \mu y\partial f(\alpha)/\partial a - n/\alpha$ =0,
$n\partial\gamma_\mu(\alpha)/\partial a = y\partial f(\alpha)/\partial a - nq/p\alpha$ =0,
If we take $n = pm$: $n\partial\gamma_\mu(\alpha)/\partial a = y\partial f(\alpha)/\partial a - qm/\alpha$ =0.
We obtain the same equation with $qm$ instead of $n$ and the same asymptotic distribution of zeros with $y = qms$.
Let $H_\infty(s) = 0$ with $\alpha$, $s$ and $\kappa$ fixed. As $H_{n-1}(s_n)$ is polynomial in $s_n$, the solution $s_n$ of the equations $\partial\gamma_n(\alpha)/\partial a = 0$ and $\Im\mathrm{m}(\gamma_n) = \pi\,\kappa_n$ is uniformly continuous function of $\kappa_n$ near to $\kappa$ for $\forall\ n$>N. With a multidimensional Dirichlet's approximation, for $\forall n$>N, $\exists k$ and $\kappa_n = k/n$ : $|\kappa - \kappa_n| < 1/n^2 < 1/N^2$. ■
As in the previous lemma, for $\forall n > N$ we have $s_n$ near to $s$ for a known critical point $\alpha$.

**Proposition**

*For $\forall n > N, s_n$ verify the recurrence: $s_{n+1}\partial f(f(\alpha))/\partial a$=$s_n$ . But asymptotically, we can have solutions of the equation: $s_{n+1}\partial f(f(\alpha))/\partial a$=$\mu s_n$ with $\mu \in \mathbb{R}$.*

*If the matrix $\partial f(f(\alpha))/\partial a$ is diagonalizable, we search the conditions on $s_n$ to be an eigen vector of this matrix. If we note eigen vector $\mu = \rho e^{i\vartheta}$, it induces:*

$$\vartheta = k\pi/n\,,$$

*and asymptotically $\vartheta = \omega\pi$ with uniform $\omega \in [0,1]$ under the condition that the vector $\kappa$ is uniform on $[0,1]^{\ r}$.*

■ At time $n$, we have: $\partial\gamma(\alpha)/\partial a = s_n\partial f(\alpha)/\partial a - 1/\alpha$ =0 ) and $H_{n-1}(s_n) = 0$. We iterate $f$ : $s_{n+1}\partial f \circ f(\alpha)/\partial a - 1/\alpha$ =0.
As: $s_n\partial f(\alpha)/\partial a - 1/\alpha$ =0,
we obtain: $s_{n+1}\partial f(f)/\partial f|_{f=f\ (\alpha)}$=$s_n$
The complex matrix $\partial f(f(\alpha))/\partial f$ doesn't depend on $n$. We can search the asymptotic behavior of $s_n$. That means:

$$s\partial f(f(\alpha))/\partial a\text{=}s.$$

There is no solution, except if an eigen value is 1.
But we have seen in the lemma that we can modify asymptotically the equation for $\forall\ \mu \in \mathbb{R}$:

$$s\partial f(f(\alpha))/\partial a\text{=}\mu s.$$

If we write $\mu = \rho e^{i\vartheta}$, can approximate the real part of $\mu$, but the imaginary part must be null: $\Im\mathrm{m}(\lambda) = \rho\sin\vartheta = 0$. If such $s$ exist, for all $s_n$ which converges to $s$, the asymptotic relation between $s_n$ and $s_{n+1}$: $s_{n+1}T = \mu s_n T = \mu^n s_1 T$ is verified. All the $\mu^n$ must be real and $\sin n\vartheta = 0$. So: $\vartheta = k\pi/n$ and asymptotically $\vartheta = \omega\pi$ with uniform $\omega \in [0,1]$ under the condition that the vector $\kappa$ is uniform on $[0,1]^r$. This may explain the self-similarity phenomenon because at each iteration, the curve is reciprocal image of:

$$k/n < \frac{(k+1)}{n+1} < (k+1)/n. \blacksquare$$

**3- Example: the Henon's curves**

Let the iteration: $f(a,b) = (a_1, b_1) = (g(a) + b, \nu a)$

with $g(a) = \delta a - \sigma a^2$,

explicitly: $a_1 = g(a) + b$, $\quad b_1 = \nu a$,

We can put: $\delta = \lambda + \lambda'$ and $\nu = -\lambda\lambda'$.

Where $\lambda$ and $\lambda'$ are the eigen value of the linear part of $f$ at 0.

We observe:

$$H_n(x,y) = \partial^n e^{xg(a)+y\nu a}/\partial a^n|_{a=0}\,\partial^n e^{xb}/\partial b^n|_{b=0},$$

$$H_n(x,y) = x^n \partial^n e^{xg(a)+y\nu a}/\partial a^n|_{a=0}.$$

The hessian of $f$ is degenerated and we have only to study $\partial^n e^{xg(a)-y\nu a}/\partial a^n|_{a=0}$.

After normalization ($x, y$ divided by $n$), we put: $x\sigma a^2 = u^2/2$ and $(x\mu + y\nu)a = \beta u$.

Then: $\beta = (x\mu + y\nu)/\sqrt{2x\sigma}$.

We recognize in $\gamma(a) = xg(a) + y\nu a - \ln a$, the PR function of the logistic map: $f(u) = \beta u - u^2/2$. We have the density of probability of presence under the conditions if the discriminant of the PR function: $\Delta = (s\lambda)^2 - 4s < 0$ and $s > 0$. This density is:

$$p(a) = \beta/(2\pi\sqrt{4a - a^2\beta^2}).$$

With: $\beta = (x\delta + y\nu)/\sqrt{2x\sigma}$

We have a piece of stochastic parabola.

But now, we consider $f \circ f = f(a_1, b_1) = (g(a_1) + b_1, \nu a_1)$,

$$f \circ f(a,b) = (g(g(a) + b) + \nu a, \nu(g(a) + b)).$$

The PR function becomes:

$$\gamma_1(\alpha_1) = x(g(g(a) + b) + \nu a) + y\nu(g(a) + b) - \ln a - \ln b.$$

Its critical point is: $\partial\gamma_1(\alpha_1)/\partial a = \{x\, g'(g(a) + b) + y\nu\}g'(a) - 1/a = 0$.

$$\partial\gamma_1(\alpha_1)/\partial b = \{x\, g'(g(a) + b) + y\nu\} - 1/b = 0.$$

Then, fixing $x_1 = x\, g'(g(a) + b) + y\nu$, we have yet: $x_1 g'(a) - 1/a = 0$ and $x_1 - 1/b = 0$; which determines the critical point $\alpha = (a, b)$ of $\gamma(a)$ for $f$ where $b$ is unknown and $a$ is invariant. But, with $\partial\gamma_1(\alpha_1)/\partial b = 0$, we have a new equation of $b$ function of $a$:

$$bx\, g'(g(a) + b) + \nu y b + 1 = 0.$$

As: $g'(a) = \delta - 2\sigma a$, $g'(g(a) + b) = \delta - 2\sigma(g(a) + b)$; the equation in $b$ is:

$$x\, b\big(g'\big(g(a)\big) - 2\sigma b\big) + y\nu b + 1 = 0.$$

$$-2\sigma x b^2 + (x\,\delta + y\nu - 2\sigma x(\mathcal{R}e(g(a) + i\Im m(g(a)))b + 1 = 0.$$

But:

$$\partial f(f(\alpha))/\partial a = \begin{bmatrix} \delta - 2\sigma g(a) & 1 \\ \nu & 0 \end{bmatrix},$$

and the eigen vector $s = (\mu, \nu)$ verifies: $\mu\big(\mu - \delta.2\sigma g(a)\big) - \nu = 0$.

The sequence $s_n = (\mu^n, \nu\, \mu^{n-1})$ If we write $\mu = \rho\, e^{-i\theta}$, the imaginary part of $\mu$ will be null. That means: $n\theta = k\pi/n \to \omega$ uniform on (0,1). But, before the uniformity, the self-similarity dominates. $\blacksquare$

## E- The complex iterations

**1- The trinomial iterations**

We study here only trinomial iterations as the following.

Let the application of a set $C \subset \mathbb{C}$ in $C$ :

$$z_1 = f(z) = \alpha z + z^m/m;\ m > 1$$

$z, \alpha \in \mathbb{C}$ ,and $f$ is trinomial. All the following calculus are formal.

The fixed points are: $z_0 = 0$ et $z_{m0} = ((1-\alpha)m)^{1/(m-1)}$. As $\alpha \in \mathbb{C}$, we note $\alpha = re^{i\sigma}$ and $z = \rho e^{i\vartheta}$ and $\vartheta+\sigma=\omega$. We have $m-1$ fixed points $z_{m0}$:

$$z_{m0k}{=}z_{m0}e^{i2k\pi/(m-1)}.$$

The derivatives of $f$ are $\lambda_0 = \alpha$ at point 0 and $\lambda_{m0} = \alpha + m(1-\alpha)$ at points $z_{m0}$ .

The steepest descent method is supposed yet valid. We observe a great variety of behaviors.

**- *Study at 0***

The Plancherel-Rotach's function is at 0 : $\gamma(z) = sf(z) - \ln(z) = s(\alpha z + z^m/m) - \ln(z)$, the critical points are given by : $\partial\gamma(z)/\partial z = s\partial f/\partial z - 1/z = s(\alpha + z^{m-1}) - 1/z = 0.$

$$\partial\gamma(z)/\partial z = z^m + \alpha z - 1/s = 0.$$

The discriminant of this equation is: $\Delta = (-1)^{m-1}((m-1)^{m-1}(\alpha)^m + (m)^m\ (1/s)^{m-1})$.

If $\alpha$ is real and $\Delta$ is negative, or if $\alpha$ is complex, we have a complex root. In polar coordinates:

$$\partial\gamma(z)/\partial z = s\rho^m e^{im\vartheta} + sr\rho e^{i\omega} - 1 = 0,$$

and: $\Im\mathrm{m}(\partial\gamma/\partial z) = \Im\mathrm{m}(z^m + \alpha z) = 0,$

$$\rho^{m-1}\sin m\vartheta + r\sin\omega = 0.$$

We have $\rho$ function of $\vartheta$: $\rho^{m-1} = -r\sin\omega/\sin m\vartheta > 0$ with condition: $|\rho^{m-1}\sin m\vartheta| \leq r$.

But we have also: $\mathcal{Re}\left(s\rho^m e^{im\vartheta} + sr\rho e^{i\omega}\right) = 1,$

$$sr\sin(m\vartheta - \omega) = \sin m\vartheta\ .$$

We can write: $s = \sin m\vartheta\ /\ r\sin(m\vartheta - \omega) = \varphi(\theta)/r.$

Let $\chi$ uniform on (0,1). If $s$ is real, the distribution of probability of $s$ is:

$$\pi d\chi = \Im\mathrm{m}(f(z))ds;$$

As: $\Im\mathrm{m}(\partial\gamma/\partial z) = \Im\mathrm{m}(z^m + \alpha z) = 0,$

$$\pi d\chi = (m-1)/m)\ \Im\mathrm{m}(\alpha z)\ ds.$$

Then we have: $\pi d\chi = (m-1)/m) r\rho\sin\omega ds,$

and we can express $\chi$ in function of $\vartheta$ and get the low of $\vartheta$:

$$\pi d\chi = (m-1)/m)(-r\sin\omega/\sin m\vartheta)^{1/(m-1)}\sin\omega\ \ \varphi'(\theta)d\theta.$$

If, now we iterate $k$ times $f$, we have to study the asymptotical equation:

$$s\mathcal{Re}(\partial f(f)/\partial f|_{f=f(\alpha)}) = s = s\mathcal{Re}\ ((\alpha z + z^m/m)^{m-1} + \alpha),$$

where z is critical point of $\gamma$. Let $\mu = \mathcal{Re}\ ((\alpha z + z^m/m)^{m-1} + \alpha)$.

The separation between the bounded iteration and the divergence is given by:

$$|\mu| = 1$$

But, we must have also the similarity's condition: $\Im\mathrm{m}(\mu^k) = 0$.

**- Study at the points $z_{m0}$**

We have $m-1$ fixed points $z_{m0}$: $z_{m0k}{=}z_{m0}e^{i2k\pi/(m-1)}$ .

It suffices to translate the calculus at these points. The only novelty is the domain of attraction of each fixed point: (0, $z_{m0k}$). The relation of domination is given by comparison of the real part of the Plancherel-Rotach's functions.

**2- The Fatou-Julia's sets**

**-*m*=2,** if $0 > \alpha \in \mathbb{R}$**,** we recognize the logistic map when $4s + (\alpha s)^2 < 0$. If $\alpha \in \mathbb{C}$**,** we note: $\alpha = re^{i\sigma}$ and $z = \rho e^{i\vartheta}$. The critical points are solution of: $\partial\gamma(z)/\partial z = s(\alpha z + z^2) - 1 = 0.$

The imaginary part is: $\rho\sin(2\vartheta) = -r\sin(\sigma+\vartheta).$

The real part is: $s\ \rho\ (r\cos(\sigma+\vartheta) + \rho\cos(2\vartheta)) = 1.$

And: $s\ \rho r\ (\sin(2\vartheta)\cos(\sigma+\vartheta) - \sin(\sigma+\vartheta)\cos(2\vartheta)) = \sin(2\vartheta)$

$$s = -\sin(2\vartheta)^2/r^2\sin(\vartheta+\sigma)\sin(\vartheta-\sigma)$$

We write $s = (-1/r^2)\varphi(\vartheta,\sigma).$

The distribution of probability of $s$ is:

$$\pi d\chi = \Im m\ (\alpha z/2)\ ds = -r^2 \sin(\sigma + \vartheta)^2 / 2sin\,(2\vartheta)\ ds,$$

$$\pi d\chi = \sin(\sigma + \vartheta)^2 \varphi'(\vartheta, \sigma) d\vartheta / 2sin(2\vartheta).$$

It is more difficult than the real case where $\alpha \in R$.
If, now we iterate $k$ times $f$, we have to study the asymptotical equation:

$$s_n = s_{n+1} \mathcal{R}e\ (\alpha + \alpha z + z^2/2),$$

where z is critical point of $\gamma$. Let $\mu = \mathcal{R}e\ (\alpha + \alpha z + z^2/2)$.
The separation between the bounded iteration and the divergence is given by:

$$|\mu| = 1.$$

But, we must have also the similarity's condition: $\Im m(\mu^k)=0$,

$$\Im m\ (\mu) = \Im m\ (\alpha + \alpha z/2) = r(sin(\sigma) - r\ sin\,(\sigma + \vartheta)/2\ sin(2\vartheta)).$$

- $\boldsymbol{m} = \mathbf{3}$,

In this case, we have analogues of the derivatives of the Airy's function or the Hopf's bifurcation the differential iteration: $f(z)=z(1+\delta(\alpha+b|z|^2)$, $z\in \mathbb{C}$ where δ is the path.
We study: $f(z) = \alpha z + z^3/3$. The critical point of $\gamma(z) = s(\alpha z + z^3/3) - \ln(z)$, is solution of: $\partial\gamma/\partial z = z^3 + \alpha z - 1/s = 0$. If $0 > \alpha \in \mathbb{R}$;
the two imaginary solutions are obtained with the negative discriminant: $\Delta = 4\alpha^3 + 27/s^2 < 0$.

## Part II: Differential ordinary equations as iterations

### A- The probabilistic approach of the Bendixon- Poincaré problem

In this part, we apply the results on the bounded iteration to the ordinary differential equation. Partial differential equations are more difficult to study and it will be the subject of a next paper. Here, we can say, under good conditions, that a bounded ODE is deterministic near the origin of the process, but may have stochastic or fixed cycles after a very long time.
We correct a previous paper with new demonstrations. [4]
With the probabilistic method, we obtain some new results in mechanic or in physics, but we meet also many new difficulties due to the particular steepest descent's method used to study the Plancherel-Rotach's function.
We consider the ordinary differential equation ODE:

$$da/dt = F(a),$$

where $a \in C \subset \mathbb{R}^d$ or $\mathbb{C}^d$, $t \in \mathbb{R}^+$, $F(a)$ is a $C^\infty$-application of $a \in C$ in C. The domain C is supposed bounded. The problem is to find a function $a(t)$ verifying this equation with an initial condition: $a(t_0) = a_0$.

#### 1- The differential iteration and the asymptotic invariant measure

We associate to the ODE, the differential iteration $f(a)$ belonging in a bounded domain C:

$$f(a) = a + \delta F(a)$$

$$a_1 = a + \delta F(a)$$

where $\delta = t/n$ is the path. We recall that the fixed points of a differential iteration are the zeros $\alpha$ of $F$: $F(\alpha) = 0$.

**Lemma 4**

*The invariant measure P of the differential iteration verifies:*

$$\int e^{ya}(yF(a))^n dP(a) = 0, \forall n.$$

*So, for* $\forall\tau > 0$: $E(e^{yA}) = E\left(e^{y(A+\tau F(A))}\right)$.

*The invariant measure Pis valid for all* $\tau > 0$ *instead of the initial small* δ.

■ For $\delta \neq 0$, we write the Fourier-Laplace's transform of the invariant measure $P_\delta$ associated to the differential iteration $a + \delta F(a)$ :

$$E(e^{yA}) = E\left(e^{y(A+\delta f(A))}\right).$$

And: $\theta(y,\delta)$=$E\left(e^{yA}(1-e^{y\delta F(A)}\right)$ =0.

When $\delta \to 0$: $\theta(y,\delta)$=0=$E\left(e^{yA}(1-e^{y\delta F(A)}\right) \approx \delta E(yF(A)e^{yA})$,

$$yE(F(A)e^{yA}) = y\int e^{ya}F(a)dP_\delta(a) = 0.$$

By derivation on $\delta$, we have also:

$$\int e^{ya}(yF(a))^n dP(a) = 0\ ,\forall n.$$

So, $P_\delta$ doesn't depend on $\delta$ near by $\delta$ =0. and $P_\delta \to P$ . Let $\Delta = supp(P)$ its support.

Then, by summation with arbitrary $\tau^n/n!, we$ , we have the result. ■

This lemma is consistent with the definition of the differential equation $\mathrm{d}a/\mathrm{d}t = F(a)$ when we put: $t = ct'$ ; the equation becomes $\mathrm{d}a/\mathrm{d}t' = \mathrm{c}F(a)$ where $c$ is as large as we want.

Under our hypothesis, the theoretical Caratheodory's solution exists $a(t)$ for some $t > t_0$ with $a_0 = a(t_0)$ and $a_0$ belongs to a neighborhood of the initial point.

$$a(t)= a_0+\int_{t_0}^{t} F\big(a(u)\big)\mathrm{d}u.$$

When we iterate $n$ times $f$, we have the iteration:

$$a_n = f^{(n)}(a(t_0)).$$

It can be written: $a_n = f^{(n)}(a_0) = a_0 + \delta\left(\Sigma_{p=0}^{p=n-1}F(a_p)\right)$

When $n \to \infty$ , this solution $a_n \to a(t)$ converges to the Caratheodory's solution:

$$\lim_{n\to\infty} f^{(n)}(a_0)=a(t) = a_0+\int_{t_0}^{t} F\big(a(u)\big)\mathrm{d}u$$

A priori, this deterministic solution depends on $a_0$. In this case, we have stopped the iteration at the step $n$ . But, if we forget to stop the computer at this step, as does Lorenz, the sequence of the iterates $a_p$ increases indefinitely with $p$ , for the constant step $\delta$. Now, the time $t$ doesn't have any sense. We obtain a good approximation of the behavior of $a(\tau)$ when $\tau = p\delta \to \infty$ with $p$ . So, we can identify two behaviors for a differential equation: the first, deterministic for a finite time, the second potentially probabilistic for an infinite time. The behavior at the infinite may be the extension of the Caratheodory's solution when $t \to \infty$ ; and this is the deterministic point of view which is relevant sometimes. But, we can see some curious erratic behaviors which are not deterministic. We suggest to study it with the Fourier-Laplace's transform of the invariant measure $P$ associated to the differential iteration.

**Lemma 5**

*Let $a(0)$ a point of the support of the invariant measure and $a(t) \neq a(0)$ the unique solution of Carathéodory for $t \neq 0$ . We have:* $E\left(e^{y(A(0))}\right) = E\left(e^{yA(t)}\right)$.

*It exists an almost period $\tau$ such as $A(0) = A(\tau)$.*

■ If we iterate $n$ times the differential iteration from the initial position $a(0)$ at $t = 0$ :

$$a(\tau) = a(0)+\sum_1^n F(a(k\tau/n))\tau/n.$$

As $f$ is continue and bounded, the limit of this sum is the integral when $n \to \infty$:

$$\lim_{n\to\infty} \sum_1^n F(a(k\tau/n))\tau/n = \int_0^\tau F(a(t))dt\ .$$

So: $a(\tau) = a(0) + \int_0^\tau F(a(t))dt.$

If $a(0) \in \Delta$ support of $P$, the definition de PF for the $n$ −th iterate $a(n)$:

$$E\left(e^{y(A(0))}\right) = E\left(e^{yA(n)}\right).$$

Then, the theorem of the dominated convergence implies:

$$E\left(e^{yA(\tau)}\right)=\lim_{n\to\infty} E\left(e^{yA(n)}\right).$$

And: $E\left(e^{y(A(0))}\right) = E\left(e^{y(A(0))+\int_0^{\tau} F(A(t))dt)}\right)$,

$$E\left(e^{y(A(0))}(1 - e^{y(\int_0^{\tau} F(A(t))dt)}\right) = 0.$$

The domain $C$ is bounded, and we can translate $C$ so that all the coordinates of $C$ are positive and $F(C)$. Taking $y$ positive, the expression is always negative; so, almost everywhere for the invariant measure:

$$1 - e^{y(\int_0^{\tau} F(a(t))dt)} = 0,$$

and: $y\int_0^{\tau} F(a(t))dt = 0$ . ■

Now, we have to define what are the values of $\tau$ when this can happen.

## 2- The Fredholm's alternative[3] & [4]

**Proposition**

*When the number of iterations $n \to \infty$ , with $s = y/n$ the critical point is defined by:*

$$s + \tau s\partial F(\alpha)/\partial a - 1/\alpha = 0$$

*If $s_\alpha$ is a particular solution and if $\vec{s}$ is an eigenvector of $\partial F(\alpha)/\partial a$ for the eigenvalue* $-1/\tau$, *the general solution is disjunctive:*

$$s = s_\alpha \text{ if } \tau \neq -1/\lambda_\alpha \text{ or } s = \vec{s} \text{ if } \tau \neq -1/\lambda_\alpha.$$

*The eigenvalue* $-1/\tau$ *can be interpreted as a critical asymptotic frequency.*

*If this period exists, it may be at stochastic as $\alpha$.*

*If such a period doesn't exist, then the process flattens in the kernel of $\partial F(\alpha)/\partial a$ .*

■ The asymptotic distribution is defined by the limit of $e^n(y) = \partial^n(e^{ya} - e^{yf(a)})/\partial a^n$ when $n \to \infty$, for $a = 0$. the Cauchy's representation gives:

$$e^n(y) = c\oint_\Gamma \frac{e^{ya} - e^{yf(a)}}{a^n} da = c\oint_\Gamma e^{ya - n\ln a}(1 - e^{\delta yF(a)})da,$$

But, with the mean's formula: $1 - e^{\delta yF(a)} = -\int_0^{\delta} yF(a)e^{tyF(a)}dt$.

As we only need the zeros of $e^n(y)$, the research is limited to find the zeros of:

$$\oint_\Gamma yF(a)e^{ya + tyF(a) - n\ln a}\, da = 0,$$

With $t$ arbitrary as δ. As δ is "elastic" with the lemma 4, we note δ instead of $t$.

Let a point $a(\tau_0) \in \Delta$ support of $P$; as $a(\tau_0)$ has a critical point $\alpha_0$ defined by at time $\tau_0 = n\delta$:

Putting $y = ns_0$, the Plancherel-Rotach function is:

$$n\gamma(a(\tau_0)) = ns_0\big(a + \delta F(a)\big) - n\ln(a).$$

$$\partial\gamma(\alpha_0, \tau_0)/\partial a = s_0 + \delta s_0 \partial F(\alpha_0)/\partial a - 1/\alpha_0 = 0.$$

$$\Im m(\gamma(\alpha_0, \tau_0) - \theta_0 = \chi_0 \pi.$$

Now, we start from $a(\tau_0)$ and we iterate $k$ times more $f(a)$ with $\tau = k\,\delta$.

The Plancherel Rotach function for the $(n + k)th$ iterated $a((n + k)\delta) = a(\tau_0 + \tau)$ is now:

$$(n + k)\gamma(a(\tau_0 + f(a) -) = (n + k)s'(a + \delta F(a)) - (n + k)\ln(a).$$

Then, $\gamma(a(\tau_0 + \tau)) = \gamma(a(\tau_0))$ under the only change: $ns_0 = (n + k)s'$.

At the opposite of the previous study, we don't write the critical point $\alpha$ as a function of $s$, but $s$ as a function of $\alpha$.

**Recall**: *For fixed $\alpha$, we recognize in the equation $\partial\gamma(\alpha)/\partial a = 0$ the linear affine equation of $s$ depending on the parameter $\delta$:*

$$s + \delta s\partial F(\alpha)/\partial a - 1/\alpha = 0 .$$

*In this case, the solution of this equation of Fredholm's type is alternative:*

*- if $s + \delta\partial f(\alpha)/\partial as = 0$ has no solution other than $\delta = 0$ , then $s = 1/a$ ;*

*- if not: $s + \delta\partial f(\alpha)/\partial as = 0$ has solutions where $s(\alpha)$ is eigen vector of $\partial f(\alpha)/\partial a$ for each eigen value $\lambda$ de $\partial f(\alpha)/\partial a$ at $\alpha$.*

*Then, if we have a particular solution $s_a$ :*

$$s_\alpha + \delta s_\alpha \partial F(\alpha)/\partial a - 1/\alpha = 0.$$

*Formally:* $\quad s_\alpha = (Id + \delta\partial F(\alpha)/\partial a)^{-1}(1/\alpha).$

*This solution $s_\alpha$ is valid for all $\delta \neq$-1/$\lambda_\alpha$ where $\lambda_\alpha$ is eigenvalue of $\partial F(\alpha)/\partial a$ at $\alpha$. The general solution will be $s$ =$s_\alpha$ +$\vec{s}$ , with: $\vec{s} = -\delta\vec{s}\partial F(\alpha)/\partial a$. As $\delta \neq$-1/$\lambda_\alpha$, $\vec{s}$ =0. But, as δ is increasing from 0 to ∞, the general solution is disjunctive and shows discontinuities at the eigenvalues $\lambda_\alpha$. More $s_\alpha$ is orthogonal to the eigen vectors of $\partial F*(\alpha)/\partial a$.*

Here, the critical point $\alpha_0$ is the same for $\gamma(a(\tau_0)$ and for $\gamma(a(\tau_0 + \tau)$. So:

$$s_0 + \delta s_0 \partial F(\alpha_0)/\partial a - 1/\alpha_0 = 0,$$

and: $\quad s' + \delta s' \partial F(\alpha_0)/\partial a - 1/\alpha_0 = 0.$

By difference: $\quad (s' - s_0) + \delta(s' - s_0)\partial F(\alpha_0)/\partial a = 0.$

As $(n + k)s' = ns_0$, we have: $s' + \delta s' \partial F(\alpha_0)/\partial a$=0.

The only possible solutions are eigen vectors for eigen values $\lambda$ of $\partial F(\alpha_0)/\partial a$ and $\delta$ =*-1* /$\lambda$. Then, as δ is positive, the only possible periods are defined by the eigen values $\lambda$ real and negative but each eigen vector must be orthogonal to (1/$\alpha_0$). If it doesn't exist, $s'$ is eigen vector of: $\lambda s' - s'\partial F(\alpha_0)/\partial a$=0; but $s' + \delta s' \partial F(\alpha_0)/\partial a$=0 and if δ → ∞, then $s'\partial F(\alpha_0)/\partial a$=0; $s'$ belongs in the kernel of $\partial F(\alpha_0)/\partial a$. In general, $s' = 0$. The process is flattening and change. ■

**Remarks**

- When we have many fixed points, the complete solution is more difficult to apprehend because we meet some problems with the transitions from a domain of a fixed point to an another when all the domains are repulsive (as in the Lorenz's equations).
- Assuming the repulsive case, in the §C we study only the zeros of $H_n(y)$ instead of $e^n(y) = \partial^n(e^{ya} - e^{yf(a)})/\partial a^n|_{a=0}0$,. As we observe $e^{ya} - e^{yf(a)} = -\int_0^1 y(f(a) - a)e^{ya+ty(f(a)-a)}dt$, and the integrand $e^{ya+ty(f(a)-a)}$ is the same as for a differential iteration with $F(a) = f(a) - a$.

That means for an arbitrary bounded iteration $f(a)$, we may have stochastic periodic behavior with period $\tau = 1/(1 - \lambda)$, where $\lambda$ is eigen value of $\partial f(\alpha)/\partial a$ or a flattening process.

### 3- The interaction between the asymptotic periodic cycles

This study corresponds to the asymptotic analysis of $f^{(p)}$, ∀ $p$.

We have seen the recurrence of $s_n$ at the critical point $\alpha$: $s_{n+1}\partial f(f(\alpha))/\partial a$=$\mu s_n$ *with* $\mu \in \mathbb{R}$.

Here, the differential iteration is: $f(a) = a + \delta F(a)$

We obtain: $\quad s_{n+1}(I + \delta\partial F(f(\alpha))/\partial a)$=$\mu s_n$

As δ → 0: $\quad s_{n+1}(I + \delta\,\partial F(\alpha)/\partial a)$=$\mu s_n$

And we have to search an asymptotic $s$ eigen vector of:

$$s(I + \delta\,\partial F(\alpha)/\partial a) = \mu s$$

$$s\partial F(\alpha)/\partial a = ((\mu - 1)/\delta)s = \lambda s$$

## B- Application to the mechanics

**Example 1**: ***The Hamilton's equations*** *in* $\mathbb{R}^d$ *for an autonomous system.*

Let the equations: $\quad \mathrm{d}p/\mathrm{d}t = -\partial H(p,q)/\partial q,$

$\quad \mathrm{d}q/\mathrm{d}t = \partial H(p,q)/\partial p$ , $\qquad a = (p,q) \in \mathbb{R}^{2d}$

we note: $\quad F(a) = (-\partial H(a)/\partial q, \partial H(a)/\partial p),$

then: $\quad da/dt = F(a),$

the corresponding iteration: $f(a) = a + \delta F(a)$, with $\delta = t/n$.
where $a \in C \subset \mathbb{R}^{2d}$, $t \in \mathbb{R}^{+}$, $F(a)$ is a $C^{\infty}$-application of $a \in C$ in C. The domain C is supposed bounded. We have to search the fixed points defined by:

$$\partial H(p,q)/\partial q = 0,$$
$$\partial H(p,q)/\partial p = 0.$$

We study now the asymptotic probabilistic solution of the differential iteration with the PR integrand:

$$\gamma(a) = xp + yq - x\delta\, \partial H(p,q)/\partial q + y\delta \partial H(p,q)/\partial p - \ln(p) - \ln(q).$$

The critical point is for $\delta = t$:

$$\partial\gamma(a)/\partial p = x - xt\, \partial^2 H/\partial q\partial p + yt\partial^2 H/\partial p^2 - 1/p = 0,$$
$$\partial\gamma(a)/\partial q = y - xt\, \partial^2 H/\partial q^2 + yt\partial^2 H)/\partial p\partial q - 1/q = 0.$$

Let the Hamiltonian matrix:

$$\partial^2 \widetilde{H} = \begin{pmatrix} -\partial^2 H/\partial q\partial p & \partial^2 H/\partial p^2 \\ -\partial^2 H/\partial q^2 & \partial^2 H)/\partial p\partial q \end{pmatrix}$$

With the antisymmetric matrix: $J = \begin{bmatrix} 0 & 1 \\ -1 & 0 \end{bmatrix}$, $\partial^2 \widetilde{H} J$ is symmetric.
If $1/a$=$(1/p, 1/q)$ and $s = (x, y)$ :

$$st\partial^2 \widetilde{H} + s = 1/a$$

If we have a particular solution resolving this equation, it remains a linear asymptotic system:

$$st\partial^2 \widetilde{H} + s = 0.$$

The eigen vector corresponding to the smallest eigen value when $t \to \infty$, gives $s\partial^2 \widetilde{H} = 0$.
Asymptotically, the projection of the Laplacian on this eigen vector is null:

$$\alpha\partial^2 H/\partial p^2 + \beta\partial^2 H/\partial q\partial p + \gamma\partial^2 H/\partial q^2 = 0$$

**Example 2:** ***In the case of a potential $V(q)$ in $\mathbb{R}^d$ with derivatives decreasing at $\infty$:***

$$H(p;q) = V(q) + p^2/2m$$

The Hamilton's equations are: $\mathrm{d}p/\mathrm{d}t = -\,\partial V(q)/\partial q$,
$\mathrm{d}q/\mathrm{d}t = p/m$.
The asymptotic Plancherel-Rotach's function is:

$$\gamma(p,q) = xp + yq - x\delta\partial V(q)/\partial q + \delta yp/m - \ln(p) - \ln(q)$$

We have two independent equations $\gamma(p,q) = \gamma_0(p) + \gamma(q)$
with: $\gamma_0(p) = xp + \delta yp/m - \ln(p)$,
and: $\gamma(q) = yq - y\delta\partial V(q)/\partial q - \ln(q)$.
The critical point is: $x + \delta y/m - 1/p = 0$ ,
and: $y - x\delta\partial^2 V(q)/\partial q^2 - 1/q = 0$.
The fixed point is defined by $p = 0$ and $\partial V(q)/\partial q = 0$.
The solution depends only on $V(q)$ with independency of $\gamma(q)$ and $\gamma_0(p)$.
If we have a particular solution resolving this equation, it remains a linear system:

$$x + \delta y/m = 0 \qquad \text{and} \qquad y - x\delta\partial^2 V(q)/\partial q^2 = 0.$$

So: $y - m\delta^2 y\partial^2 V(q)/\partial q^2 = 0$.
It means: $y$ is eigen vector of the real symmetric matrix $\partial^2 V(q)/\partial q^2$ . As $t \to \infty$, if the phenomenon is bounded and cyclical $\delta \to \infty$: $y\partial^2 V(q)/\partial q^2 = my/\delta^2 \to 0$.
So, we must take in priority the eigen vectors corresponding to the smallest eigen values:

$$y\mathrm{d}p/\mathrm{d}t = -\, y\partial V(q)/\partial q = \mathrm{C}.$$

In $\mathbb{R}^d$, the projection of the forces on the eigen vectors $y$ is constant; taking $q \to \infty$ , C=0. Along the vector $y$, the behavior of the bodies is linear as the barycenter. Asymptotically the cyclical behavior exists only in few dimensions corresponding to the smallest eigen values. We can imagine for very small bodies a unique vibration along only each eigen vector.

**Example 3:** ***The n-bodies problem: the solar system***

We have a sun with a mass $M$ and $n$-1 planets with a total mass $\mu = M/1500$. The heaviest planet of the system is Jupiter with a mass $M/1000$. So, we call $\varepsilon_\ell$ the mass of the planets. Let $a$ be a point of $\mathbb{R}^d$; $a_1$ and $a_\ell$ are the positions of the sun and the planets; $a(i,\ell)$ is the distance between two bodies; the potential function is:

$$V(q) = M\Sigma_{\ell=2}^{\ell=n}\, \varepsilon_\ell g(a(1,\ell)) + \Sigma_{\ell>1}^{\ell=n}\Sigma_{i>1}^{i=n}\varepsilon_\ell\varepsilon_i g(a(i,\ell)),$$

with $|a_\ell| = \|a_\ell - a_1\|^2 = a(1,\ell)\,;\, g(|a_\ell|) = 1/|a_\ell|$ and $q = a$. (But $g$ can be more general). The main difficulty is the discontinuity at $|a_\ell| = 0$. We suppose the validity of the calculus and weak inter actions: then, all the second member $\Sigma_{\ell>1}^{\ell=n}\Sigma_{i>1}^{i=n}\varepsilon_\ell\varepsilon_i g(a(i,\ell))$ is negligible and:

$$V(q) = M\Sigma_{\ell=2}^{\ell=n}\, \varepsilon_\ell g(|a_\ell|).$$

The fixed points of the iteration are defined by: $\partial V(a)/\partial a = 0\,, p = 0$ where the forces and the celerity are null. They are the Lagrange's points.

Let: $\partial V(a)/\partial a = M\Sigma_{\ell=2}^{\ell=n}\, \varepsilon_\ell \partial g(|a_\ell|)/\partial|a_\ell|\; \partial|a_\ell|/\partial a.$

That means: $\partial|a_\ell|/\partial a_1 = 2(a_1 - a_\ell)$ and $\partial|a_\ell|/\partial a_\ell = 2(a_\ell - a_1) = -\partial a/\partial a_1,$

$$\partial^2|a_\ell|/\partial {a_1}^2 = \partial^2|a_\ell|/\partial {a_\ell}^2 = 2 \text{ and } \partial^2|a_\ell|/\partial a_\ell\partial a_1 = -2,$$

$$\partial V(a)/\partial a_1 = M\Sigma_{\ell=2}^{\ell=n}\, \varepsilon_\ell g'(|a_\ell|)\partial|a_\ell|/\partial a_1,$$

$$= 2M\Sigma_{\ell=2}^{\ell=n}\, \varepsilon_\ell g'(|a_\ell|)(a_1 - a_\ell).$$

$$\partial V(a)/\partial a_\ell = M\varepsilon_\ell g'(|a_\ell|)\partial|a_\ell|/\partial a_\ell.$$

$$\partial^2 V(q)/\partial {a_1}^2 = M\partial(\Sigma_{\ell=2}^{\ell=n}\, \varepsilon_\ell g'(|a_\ell|)\partial|a_\ell|/\partial a_1)/\partial a_1,$$

$$= M\Sigma_{\ell=2}^{\ell=n}\varepsilon_\ell(g''(|a_\ell|)(\partial a/\partial a_1)^2 + 2g'(|a_\ell|))$$

$$\partial^2 V(q)/\partial {a_\ell}^2 = -\partial^2 V(a)/\partial a_1\, \partial a_\ell\,,$$

$$= M\varepsilon_\ell\, (g''(|a_\ell|)(\partial a/\partial a_1)^2 + 2\, g'(|a_\ell|))\,.$$

We note: $c_\ell = \partial^2 V(q)/\partial {a_\ell}^2$ then $c_1 = \Sigma_{\ell=2}^{\ell=n} c_\ell$, and: $\partial^2 V(a)/\partial a_1\, \partial a_\ell = -c_\ell$,

$$\partial^2 V(a)/\partial a_l\, \partial a_\ell = 0 \quad \text{for} \quad l \neq 1,\ell.$$

We observe: $\partial^2 V(a)/\partial a_1\, \partial a_\ell + \partial^2 V(a)/\partial {a_\ell}^2 = 0.$

By integration: $\partial V(q)/\partial q_1 + \partial V(q)/\, \partial q_\ell = C\,, \forall q_1 q_\ell\,.$

Then $C = 0$ and the problem is reduced to the two-body problem for the $n$-1 planets.

We have the asymptotic critical equation: $y - \tau y \partial^2 V(q)/\partial q^2 = 0.$

The determinant $\Delta$ of $\partial^2 V(q)/\partial q^2$ is null. What are the eigen values other than $\tau = 0$?

Let the determinant: $\Delta(\lambda) = |\partial^2 V(q)/\partial q^2 - \lambda I|.$

Let $\Pi_n(\lambda) = \Pi_{\ell=2}^{\ell=n}(\lambda - c_\ell)$, $\tau = 1/\lambda$ and $\lambda^{-n}\Pi_n(\lambda) = \tau^n\Pi_n(1/\tau)$ is a polynomial in $\tau$ with degree $n$.

**Proposition**

$$\Delta(\lambda) = (-1)^{n+1}\lambda^{2+n} d\lambda^{-n}\Pi_n(\lambda)/d\lambda$$

*So, the eigen values $\lambda_\ell$ of $\Delta(\lambda)$ give the extrema of the polynomials $\tau^n\Pi_n(1/\tau)$ and $c_\ell$ are the roots of $\Pi_n(\lambda)$. Each eigen value is between two consecutive roots: $c_\ell < \lambda_\ell < c_{\ell+1}$.*

*Asymptotically the cyclical behavior exists only in few dimensions corresponding to the smallest eigen values. All the forces and the orbits become flatter on a plane.*

■ As $\partial^2 V(q)/\partial q^2$ is symmetric real, the eigen values are real and the question is if the eigen value is more or less than 1.

We have:

$$\Delta(\lambda) = \begin{vmatrix} c_1 - \lambda & -c_2 & -c_n \\ -c_2 & c_2 - \lambda & 0 \\ \vdots & \vdots & \vdots \\ -c_n & 0 & c_n - \lambda \end{vmatrix}$$

with $c_1 = \Sigma_{\ell=2}^{\ell=n} c_\ell$.

Adding all columns in the first column, then getting $-\lambda$ in factor, we obtain by summation on the first row:

$$= (-\lambda) \begin{vmatrix} n & -\lambda & -\lambda \\ 1 & c_2 - \lambda & 0 \\ 1 & 0 & c_n - \lambda \end{vmatrix} = (-1)^n \lambda \begin{vmatrix} n & \lambda & \lambda \\ 1 & \lambda - c_2 & 0 \\ 1 & 0 & \lambda - c_n \end{vmatrix}$$

Then: $\Delta(\lambda) = (-1)^n \lambda( n\Pi_n(\lambda) - \lambda \mathrm{B}_n(\lambda)),$

and:

$$\mathrm{B}_n(\lambda) = \begin{vmatrix} 0 & \lambda & \lambda \\ 1 & \lambda - c_2 & 0 \\ 1 & 0 & \lambda - c_n \end{vmatrix}$$

We have the relation: $\mathrm{B}_n(\lambda)=(\lambda - c_n)\mathrm{B}_{n-1}(\lambda) + \Pi_{n-1}(\lambda)$ ,

$=(\lambda - c_n)\mathrm{B}_{n-1}(\lambda) + \Pi_n(\lambda)/(\lambda - c_n),$

$=(\lambda - c_n)(c_{n-1} - \lambda)\mathrm{B}_{n-2}(\lambda) + \Pi_n(\lambda)(1/(\lambda - c_n) + 1/(\lambda - c_{n-1})).$

So, by induction: $\mathrm{B}_n(\lambda)=\Pi_n(\lambda)(\Sigma_{\ell=2}^{\ell=n} 1/(\lambda - c_\ell)$ ,

and: $\Delta(\lambda) = (-1)^n \lambda^2 \Pi_n(\lambda)( n/\lambda - (\Sigma_{\ell=2}^{\ell=n} 1/(\lambda - c_\ell).$

But we recognize: $d\lambda^{-n}\Pi_n(\lambda)/d\lambda = -\lambda^{-n}\Pi_n(\lambda)( n/\lambda - (\Sigma_{\ell=2}^{\ell=n} 1/(\lambda - c_\ell)).$

$\Delta(\lambda) = (-1)^{n+1}\lambda^{2+n} d\lambda^{-n}\Pi_n(\lambda)/d\lambda$

If $\sigma_k$ is symmetric function of order of the $n - 1$ roots $c_\ell$ of $\mathit{\Pi}_n(\lambda)$ with $\sigma_0 = 1$ :

$\mathit{\Pi}_n(\lambda)=\Sigma_{k=1}^{k=n-1}(-1)^k \sigma_k \, \lambda^{n-1-k}$,

$d\lambda^{-n}\Pi_n(\lambda)/d\lambda = -\Sigma_{k=1}^{k=n-1}(-1)^k (k+1)\sigma_k \, \lambda^{-2-k}$,

$\Delta(\lambda) = (-1)^n \Sigma_{k=1}^{k=n-1}(-1)^k (k+1)\sigma_k \, \lambda^{n-k}$ .

We note: $\Delta = \lambda(\lambda - \lambda_2) \ldots (\lambda - \lambda_\ell) \ldots (\lambda - \lambda_n).$

The eigen values $\lambda_\ell$ give the extrema of the polynomial $\tau^n \Pi_n(1/\tau)$ where the $1/c_\ell$ are the roots of $\tau^n \Pi_n(1/\tau)$; then, we have: $c_\ell < \lambda_\ell < c_{\ell+1}$.

As we have seen, in $\mathbb{R}^d$, the projection of the forces on the eigen vectors $y$ is constant.

Except strong correlations between the small bodies (as moon with earth). In this case, we must rewrite the problem with: $V(q) = M\Sigma_{\ell=2}^{\ell=n} \varepsilon_\ell g(a(1,\ell)) + \Sigma_{\ell>1}^{\ell=n}\Sigma_{i>1}^{i=n}\varepsilon_\ell \varepsilon_i g(a(i,\ell))$

We don't must forget that the $c_\ell$ are randomized. We can use also the tools of the stellar statistics when we have many bodies. ■

**Example 4:** ***the hessian is degenerated***: ***the Lorenz's equation***[7].

Generally, the hessian is not definite negative. The Lorenz's equation is a particularly important example because the differential iteration can be broken down into three independent iterations with a remarkable feature: a partial linearity; an iteration with a negative hessian which induces a probabilistic solution and another with a positive hessian. It is an ideal example to clarify the previous results. However, as there is an interpenetration of the distributions related to each fixed point, the connection between the various results remains delicate. The probabilistic presentation seems to be the best: it gives the probability of presence except at the places where the domination changes; in this case, we go from a basin to an another as a ping-pong ball.

**The iteration at its repellent fixed points**

The vectors of this equation are written in bold notations: $\mathrm{d}\boldsymbol{a}/\mathrm{dt} = F(\boldsymbol{a})$ where $\boldsymbol{a} = (a, b, c)$:

$\mathrm{d}a/\mathrm{dt} = \sigma(b-a),$

$\mathrm{d}b/\mathrm{dt} = \rho a - b - ac,$

$\mathrm{d}c/\mathrm{dt} = -\beta c + ab.$

The differential equation applies a bounded set C in itself for $\delta > t > 0$ (the phenomenon is occurring between a cold sphere at -50° and hot sphere, the earth, at +15° as the terrestrial atmosphere is modelled by Lorenz).

*The differential iteration* $\boldsymbol{a}_1 = f(\boldsymbol{a})$ associated with a given path $\delta = t/n$ is:

$a_1 = a + \delta\sigma(b-a),$

$$b_1 = b + \delta(\rho a - b - ac),$$
$$c_1 = c + \delta(-\beta c + ab).$$

This iteration is quadratic, but has a linearity in $a$.

We recall the known results concerning the fixed points:

The fixed points are zeros of $F(\boldsymbol{a})=0$. If $\rho > 1$ and $\alpha = \sqrt{\beta(\rho - 1)}$ , it exists three fixed points: the point $\mathbf{0}=(0,0,0)$, and two others symmetric with respect to the axis of $c$:

$$\boldsymbol{\alpha}_+ = (\alpha, \alpha, \alpha^2/\beta) \text{ and } \boldsymbol{\alpha}_- = (-\alpha, -\alpha, \alpha^2/\beta).$$

At $\mathbf{0}$, the eigenvalue's equation $\lambda$ of the linear part is:

$$(\beta + \lambda)\,[(\sigma+\lambda)\,(1+\lambda) - \sigma\rho] = 0,$$

But, at $\boldsymbol{\alpha}_+$ or at $\boldsymbol{\alpha}_-$: $\lambda(\beta+\lambda)(1+\sigma+\lambda) - \alpha^2(2\sigma + \lambda) = 0.$

Coefficients $\beta, \sigma, \rho$ are such as these three fixed points are repellent; that means we have to study the distributions around each fixed point. We don't speak here about attractive cycles, resonances, and some particular values of the parameters, etc. It remains many points to clarify.

**Analysis of the hessian**

Projecting $f(\boldsymbol{a})$ on an axis $\boldsymbol{y} = (x, y, z)$, we write:

$$\boldsymbol{y}\, f(\boldsymbol{a}) = L(\boldsymbol{a}) + \delta Q(\boldsymbol{a}),$$

where $L(\boldsymbol{a})$ is linear for $\boldsymbol{a}$: $L(\boldsymbol{a}) = x(a + \delta\sigma(b-a)) + y(b + \delta(\rho a - b)) + zc\,(1 - \delta\beta)$,

with: $L(\boldsymbol{a}) = aL_1 + bL_2 + cL_3$ ,

$$L_1 = x(1 - \delta\sigma) + \delta\rho y\,,$$
$$L_2 = \delta\sigma x + y(1 - \delta),$$
$$L_3 = z\,(1 - \delta\beta),$$

and $Q(\boldsymbol{a})$ is quadratic: $Q(\boldsymbol{a}) = (zb - yc)\, a$.

In $\mathbb{R}^3$, $Q(\boldsymbol{a})$ is degenerated and not definite negative. First, we examine the matrix of $Q(\boldsymbol{a})$:

$$Q = \begin{bmatrix} 0 & z & -y \\ z & 0 & 0 \\ -y & 0 & 0 \end{bmatrix}$$

Let $\mu = \sqrt{y^2 + z^2}$ the positive eigenvalue of the characteristic equation of $Q$:

$$\mu(\mu^2 - y^2 - z^2) = 0.$$

The matrix of the eigenvectors $T$ is orthogonal and constant for all $\boldsymbol{a}$.

$$T = \frac{1}{\mu\sqrt{2}} \begin{bmatrix} 0 & \mu & \mu \\ y\sqrt{2} & -z & z \\ z\sqrt{2} & y & -y \end{bmatrix}$$

Corresponding to the diagonal matrix of the eigenvectors: $\Lambda = \begin{bmatrix} 0 & 0 & 0 \\ 0 & -\mu & 0 \\ 0 & 0 & \mu \end{bmatrix}$.

**Change of basis near 0**

- We calculate in the basis of eigenvectors directly with the Hermite's polynomials. As $T$ is orthogonal, the transposed $T'$ is also its inverse: $T' = T^{-1}$.

Then, the application $\boldsymbol{u} = T\,\boldsymbol{a}$ with $\boldsymbol{u} = (u, v, w)$ transforms:

$$\boldsymbol{y} f(\boldsymbol{a}) \longmapsto G(\boldsymbol{u}) = \boldsymbol{y} f(T'\boldsymbol{u}),$$
$$Q(\boldsymbol{a}) \longmapsto Q(T'\boldsymbol{u}) = \delta\mu(w^2 - v^2),$$
$$\boldsymbol{L}(\boldsymbol{a}) \longmapsto \boldsymbol{L}T'\boldsymbol{u}.$$

Now, in the basis $\boldsymbol{u}$, the function $\boldsymbol{y} f(T'\boldsymbol{u}) = G(\boldsymbol{u})$ is factorized into three independent functions:

$$G(\boldsymbol{u}) = g_1(u) + g_2(v) + g_3(w),$$

with: $g_1(u) = l_1 u;\ \ g_2(v) = l_2\, v - \delta\mu\, v^2;\ g_3(w) = l_3 w + \delta\,\mu\, w^2.$

Where: $l_1 = (\delta\sigma x + y(1-\delta) + z\,(1-\delta\beta))/\sqrt{2}$,

$$l_2 = (x - \delta\sigma x + \delta\rho y)y/\mu - (\delta\sigma x + y(1-\delta)) - z\,(1-\delta\beta))z/\mu\sqrt{2}\ ,$$

$l_3 = (x\text{- }\delta\sigma x + \delta\rho y)z/\mu + (\delta\sigma x{+}y(1\text{- }\delta))\text{-}z\,(1\text{-}\delta\beta))y/\mu\sqrt{2}$.

- We get 3 independent iterations:

. the first iteration $g_1$ is linear;

. the second $g_2$ is a sochastic iteration;

. the third $g_3$ remains positive, except if $l_3 = 0$ half the time.

- To calculate $l_1$, $l_2$ et $l_3$, we form:

$L\,(\boldsymbol{a}) = a(x\text{- }\delta\sigma x + \delta\rho y) + b(\delta\sigma x + y(1\text{- }\delta)) + zc\,(1\text{-}\delta\beta)$:

With: $L_1 = x(1\text{- }\delta\sigma) + \delta\rho y$ ; $L_2 = \delta\sigma x + y(1\text{- }\delta)$ ; $L_3 = z\,(1\text{-}\delta\beta)$,

Then: $$l\boldsymbol{u} = (l_1, l_2, l_3)\boldsymbol{u} = LT'\boldsymbol{u} = (L_1, L_2, L_3)\frac{1}{\mu\sqrt{2}}\begin{bmatrix} 0 & y\sqrt{2} & z\sqrt{2} \\ \mu & -z & y \\ \mu & z & -y \end{bmatrix}\boldsymbol{u}.$$

**-** Let the resolving gap $e^n(\boldsymbol{y}) = \partial(\partial^n(e^{yf(\boldsymbol{a})})/\partial\boldsymbol{a}^n)\partial\delta|_{a=0} = 0$ .

For $\forall t \leq \delta$, putting $\boldsymbol{a} = T'\boldsymbol{u}$, we have three separate coordinates:

$$e^n(\boldsymbol{u}) = T^n\partial(\partial^n(e^{yf(T'\boldsymbol{u})})/\partial\boldsymbol{u}^n)\partial\delta|_{\boldsymbol{u}=0} = 0\,,$$

$$\partial^n(e^{yf(T'\boldsymbol{u})})/\partial\boldsymbol{u}^n = \partial^n(e^{g_1(u)})/\partial u^n\,.\,\partial^n(e^{g_2(v)})/\partial v^n.\ \partial^n(e^{g_3(w)})/\partial w^n.$$

This gives: $\partial^n(e^{g_1(u)})/\partial u^n = l_1{}^n e^{g_1(u)}$,

$\partial^n(e^{g_2(v)})/\partial v^n = H_n(g_2(v))e^{g_2(v)}$,

$\partial^n(e^{g_3(w)})/\partial w^n = H_n(g_3(w))e^{g_3(w)}$,

And: $e^n(\boldsymbol{u}) = \partial\, l_1{}^n H_n(g_2(v))H_n(g_3(w))(e^{yf(\boldsymbol{T'u})})\partial\delta|_{u=0} = 0$.

***Proposition***

*The solution around the fixed point* **0** *consists of the intersection of the family of stochastic surfaces defined by:* $l_2/2\sqrt{\mu} \longmapsto low\ \beta(1/2,1/2)$ *with the surfaces* $\sigma x - y\text{-}z\beta = 0$ *and* $(\text{- }\sigma x + \rho y)z + (\sigma x\text{-}y{+}z\beta)\,y/\sqrt{2} = 0$.

■ With the same calculations of encodings and interchanging the derivations, we have:

$\partial l_1{}^n/\partial\delta = 0$; $\partial\, H_n(g_2(v))/\partial\delta = 0$; $\partial\, H_n(g_3(w))/\partial\delta = 0$.

We study separately the three expressions:

- First: $\partial\, l_1{}^n/\partial\delta = n\,(\partial l_1/\partial\delta)l_1{}^{n-1} = 0$.

Either $\partial l_1/\partial\delta = \sigma x - y\text{-}z\beta = 0$, or: $l_1 \sim (y+z)/\sqrt{2} = 0$.

- Second: the polynomial $H_n(g_3(w))$ when $w = 0$ is a Hermite's polynomial $H_n(x)$ where $x$ is $x = il_3/(\sqrt{2\delta\mu})$. this polynomial $i^n H_n(il_3/(\sqrt{2\delta\mu})$ is always positive half the time. In a general way:

$\partial\, H_n(x)/\partial\delta = nH_{n-1}(x)\,\partial\, x/\partial\delta = 0$. So: $d(l_3/\sqrt{2\delta\mu})/d\delta = 0$,

and: $l_3 \sim (xz\sqrt{2} + (y\text{-}z)\,y)/\mu\sqrt{2} = 0$ half the time.

- Third: in the case of $H_n(g_2(w))$, in addition to the solution $l_2 = 0$, we have to find the possible invariant distribution of $H_n(l_2/\sqrt{2\delta\mu}) = 0$.

Let the integrand of $n\gamma(w) = g_2(w) - n\ln w$.

When $\delta \to 0$, $l_2 \sim (x\,y\sqrt{2} + (y\text{-}z)\,z)/\sqrt{2}\mu$ with $\mu = \sqrt{y^2 + z^2}$.

By normalization of the coordinates $\boldsymbol{x} = (x, y, z) = \delta\,ns = (\delta\,nr, \delta\,ns, \delta\,nt)$, we obtain:

$l_2 \sim n\delta(rs\sqrt{2} + (s-t)t\,)/\,2(s^2+t^2)^{\frac{1}{2}} = n\delta l_2(\boldsymbol{s})$ .

$\delta\mu = n\delta^2(s^2+t^2)^{1/2} = n\delta^2\,\mu(\boldsymbol{s})$,

$n\gamma(v) = n(\delta l_2(\boldsymbol{s})v - \mu(\boldsymbol{s})(\delta v)^2 - \ln\delta v + \ln\delta)$.

Putting $\delta v = v$, we have: $n\gamma(v) = n(l_2(\boldsymbol{s})v - \mu(\boldsymbol{s})v^2 - \ln v)$**.**

We search the critical point: $d\gamma(v)/dv = l_2(\boldsymbol{s}) - 2\mu(\boldsymbol{s})v - 1/v = 0$.

The imaginary roots are: $v(\boldsymbol{s}) = l_2(\boldsymbol{s})/4\mu(\boldsymbol{s}) \pm i\sqrt{1/2\mu(\boldsymbol{s}) - l_2(\boldsymbol{s})^2/16\mu(\boldsymbol{s})^2}$.

Under the condition: $l_2(\boldsymbol{s})^2 < 8\mu(\boldsymbol{s})$:

$l_3 \sim (rt\sqrt{2} + (s\text{-}t)\, s)/\mu(\boldsymbol{s})\,\sqrt{2} = 0,$

Implies: $l_2(\boldsymbol{s}) = -(s-t)^2 \,/\, \sqrt{2}(s^2+t^2)^{1/2}.$

The condition becomes: $(s-t)^4 \,/\, (s^2+t^2)^{3/2} < 16.$

$l_1 = 0$ implies: $s + t = 0$, then: $s < 8$.

In any case, we observe that the conditions $l_3 = l_1 = 0$ allow us to express $r$ et $t$ depending on $s$ and we can write that the density of zeros of $s$ is now:

$q(s)\mathrm{d}s$ =Prob (1 zero between $s$, $s$+ d$s$) = $|\, \Im mf(\mathrm{v}(s))|\mathrm{d}s/\pi$.

$q(s)\mathrm{d}s = l_2(\boldsymbol{s})\sqrt{8\mu(\boldsymbol{s}) - l_2(\boldsymbol{s})^2}/8\pi\mu(\boldsymbol{s})\mathrm{d}s = \mathrm{d}\kappa.$

Then, $\kappa$ follows a uniform low on (0,1) with: $s + t = 0$ (or: $\sigma x - y\text{-}z\beta = 0$) and:

$x\, y\sqrt{2} + (y\text{-}z)\, z = 0.$

We also remark that the normalization doesn't affect the coefficients of the orthogonal matrix:

$T(x, y, z) = T(\delta nr,\ n\delta s,\ n\delta\, t) = T(r, s, t)$. ■

**Analysis near α+ and α−**

We now verify similar results the two other fixed points $\boldsymbol{\alpha}$+ and $\boldsymbol{\alpha}$−.

We search the distributions around the two other fixed points. To pass from the fixed point **0** to the fixed point **α**+ or **α**− , we have just to put in the iteration instead of $\boldsymbol{a} = (a, b, c)$: $\boldsymbol{a}'$+**α**+$=(a'+\alpha,\ b'+\alpha,\ c'+\alpha^2/\beta)$ or $\boldsymbol{a}''$+**α**−$=(a''\text{-}\alpha,\ b''\text{-}\alpha,\ c''+\alpha^2/\beta)$.

- *Calculation for* **α**+

So, for $\boldsymbol{a}'$+**α**+ $= \boldsymbol{a}$; $\boldsymbol{a}'_1 = \boldsymbol{a}'$ + **α**+ and $\boldsymbol{a}_1 = f(\boldsymbol{a})$ where $\boldsymbol{a}_1 = (a_1, b_1\ c_1)$ becomes

$\boldsymbol{a}_1 = \boldsymbol{a}'_1$+**α**+$= f(\boldsymbol{a}) = f(\boldsymbol{a}'$+**α**+$)$;

then: $\boldsymbol{a}'_1 = \boldsymbol{a}' + \delta F(\boldsymbol{a}'+\alpha+)$,

As: $F(\boldsymbol{a}) = (\sigma(b\text{-}\ a), (\rho a - b - ac), (-\beta c + ab)).$

$\boldsymbol{a}_1 = f(\boldsymbol{a})$ becomes for $\boldsymbol{a}$+**α**+:

$a'_1 = a + \delta\sigma(b\text{-}\ a) = a_1,$

$b'_1 = b + \delta(\rho a - b\text{ - }ac) + \delta(\text{-}\ \alpha c\text{-}a\alpha^2/\beta) = b_1 + \delta(\text{-}\ \alpha c\text{-}a\alpha^2/\beta),$

$c'_1 = c + \delta(\text{-}\beta c + ab) + \delta\alpha(a+b) = c_1 + \delta\alpha(a+b),$

The projection of $f(\boldsymbol{a})$ on an axis $\boldsymbol{y} = (x, y, z)$ can be written:

$\boldsymbol{y}\, f(\boldsymbol{a}') = xa_1 + yb_1 + \delta y(\text{-}\ \alpha c\text{-}a\alpha^2/\beta) + zc_1 + z\delta\alpha(a+b)$ ,

$\boldsymbol{y} f(\boldsymbol{a}') = \boldsymbol{y}\, f(\boldsymbol{a}) + \delta(a(z\alpha\text{-}y\alpha^2/\beta) + z\alpha b\text{-}\ y\alpha c)$ ,

and $Q(\boldsymbol{a})$ is invariant: $\boldsymbol{y} f(\boldsymbol{a}') = L'(\boldsymbol{a}) + \delta Q(\boldsymbol{a})$.

$L(\boldsymbol{a})$ is linear for $\boldsymbol{a}$: $L'(\boldsymbol{a}) = L(\boldsymbol{a}) + \delta(a\,(z\alpha\text{-}y\alpha^2/\beta) + z\alpha b\text{-}\ y\alpha c),$

$L'(\boldsymbol{a}) = aL'_1 + bL'_2 + cL'_3$ ,

with: $L'_1 = L_1 + \delta(z\alpha\text{-}y\alpha^2/\beta)$; $L'_2 = L_2 + \delta z\alpha$; $L'_3 = L_3 - \delta y\alpha$.

Then $T$ and $\Lambda$ remain invariant. The following is only a calculus.

We calculate $l'_1$, $l'_2$ et $l'_3$, with:

$L(\boldsymbol{a}) = a(x\text{-}\ \delta\sigma x + \delta\rho y) + b(\delta\sigma x + y(1\text{-}\ \delta)) + zc\,(1\text{-}\delta\beta)$ :

where $L_1 = x(1\text{-}\ \delta\sigma) + \delta\rho y$ ; $L_2 = \delta\sigma x + y(1\text{-}\ \delta)$ ; $L_3 = z\,(1\text{-}\delta\beta)$,

and: $\boldsymbol{l}'\mathrm{u} = (l'_1, l'_2, l'_3)\mathbf{u} = LT'\mathbf{u}$.

$$= (L_1 + \delta(z\alpha\text{-}y\alpha^2/\beta), L_2 + \delta z\alpha,\ L_3 - \delta y\alpha)\frac{1}{\mu\sqrt{2}}\begin{bmatrix} 0 & y\sqrt{2} & z\sqrt{2} \\ \mu & -z & y \\ \mu & z & -y \end{bmatrix}.$$

The results are modified; if $\boldsymbol{l} = (l_1, l_2, l_3)$ is related to 0 and $\boldsymbol{l}' = (l'_1, l'_2, l'_3)$ to $\alpha$+

$l'_1 = l_1 + \delta\alpha(z - y)/\sqrt{2},$

$l'_2 = l_2 + \delta\alpha((z\text{-}y\alpha/\beta)y\sqrt{2}\text{-z}\,(z+y)\,)/\,\mu\sqrt{2}$ ,

$l'_3 = l_3 + \delta\alpha((z\text{-}y\alpha/\beta)z\sqrt{2}\text{+y}\,(z+y)\,)/\,\mu\sqrt{2}$ .

The following calculations remain the same with these modifications.

- *Calculation for* **α**−

When $\boldsymbol{a}$ becomes $\boldsymbol{a}''$+**α**− the calculation is the same with the coordinates:

$$a''_1 = a + \delta\sigma(b - a) = a_1,$$
$$b''_1 = b + \delta(\rho a - b - ac) + \delta(\alpha c - a\alpha^2/\beta) = b_1 + \delta(\alpha c - a\alpha^2/\beta),$$
$$c''_1 = c + \delta(-\beta c + ab) - \delta\alpha(a+b) = c_1 - \delta\alpha(a+b).$$

It remains the problems of domination and frontiers between the various distributions attached at each fixed point. We have to go back to the original coordinates. And the solution gives only the probabilities of presence...

**The Rössler's attractor (particular case of the Lorenz attractor)**

Vectors are written in bold characters $\boldsymbol{a} = (a, b, c)$:

$$da/t = -(b + c),$$
$$db/dt = a + \sigma b,$$
$$dc/t = \beta + c(a - \rho).$$

The differential iteration for a fixed path $\delta = t/n$ is: $\boldsymbol{a}_1 = f(\boldsymbol{a})$ :

$$a_1 = a - \delta(b + c) \quad b_1 = b + \delta(a + \sigma b) \quad c_1 = c + \delta(\beta + c(a - \rho)).$$

The iteration is quadratic in $c\,a$, but has two linearities in $a$ and $b$.

The fixed points are zeros of $F(\boldsymbol{a}) = 0$.

If $c = \tau$ is the coordinate of the fixed point: $b = -\tau$, $a = \sigma\tau$, we have the equation:

$$\beta + \tau(\sigma\tau - \rho) = 0;$$

solutions are: $\tau = \frac{-\rho \pm \sqrt{\rho^2 - 4\sigma\beta}}{2\sigma}$.

If $\sigma \leq \beta \ll \rho$, it exists two fixed points fixed: $\boldsymbol{\tau} = (\sigma\tau, -\tau, \tau)$ and $\boldsymbol{\tau}' = (\sigma\tau', -\tau', \tau')$ .

At fixed point $\boldsymbol{\tau}$, the iteration becomes with $\alpha = \sigma\tau + \rho$:

$$\boldsymbol{a}_1 = \boldsymbol{a} + \delta F(\boldsymbol{a}):$$
$$a_1 = a - \delta(b + c),$$
$$b_1 = b + \delta(a + \sigma b),$$
$$c_1 = c + \delta(\tau a + c(a - \alpha)).$$

Let the PR function: $\gamma(\boldsymbol{a}) = \boldsymbol{y} f(\boldsymbol{a}) - n\ln(\boldsymbol{a})$ with $\boldsymbol{x} = (x, y, z)$ ,

$$\gamma(\boldsymbol{a}) = x(a - \delta(b + c)) + y(b + \delta(a + \sigma b)) + z(c + \delta(\tau a + c(a - \alpha))) - n\ln(\boldsymbol{a}),$$
$$\gamma(\boldsymbol{a}) = L(\boldsymbol{a}) + \delta Q(\boldsymbol{a})/2 - n\ln(\boldsymbol{a}),$$

with $Q(\boldsymbol{a}) = 2z\,ca$ and $L(\boldsymbol{a}) = x(a - \delta(b + c)) + y(b + \delta(a + \sigma b)) + z(c + \delta(\tau a - c\alpha))$ ,

$L(\boldsymbol{a})$ is linear in $\boldsymbol{a}$: $(x + \delta y + \delta\,\tau z, -x\delta + y(1 + \delta\sigma), -x\delta + z - z\delta\alpha)$ ,

$Q(\boldsymbol{a}) = 2z\,c\,a$ is quadratic symmetric form with matrix:

$$Q = \begin{bmatrix} 0 & 0 & z \\ 0 & 0 & 0 \\ z & 0 & 0 \end{bmatrix},$$

with the eigen vectors $\lambda = z$, $\lambda = -z$.

**Example 5*: Non-autonomous Hamiltonian system***

It is easy to verify if a Caratheodory's solution exists. If it exists, the non-autonomous Hamiltonian system $H(p, q, t)$ is depending on the time and we have a new equation:

$$dH(p, q, t)/dt = \partial H(p, q, t)/\partial t.$$

The Hamilton-Jacobi's equation with the method of the characteristics leads to the known Lie equations:

$$dx/dt = \partial H/\partial p\,,$$
$$dp/dt = -\partial H/\partial x - p\partial H/\partial z,$$
$$dz/dt = -H + p\partial H/\partial p.$$

The PR function is:

$$\gamma = u\left(x + \delta\frac{\partial H}{\partial p}\right) + v\left(p - \delta\left(\frac{\partial H}{\partial x} + p\frac{\partial H}{\partial z}\right)\right) + w\left(z - \delta\left(H - p\frac{\partial H}{\partial p}\right)\right) - \ln(p) - \ln(x) - \ln(z).$$

The critical point is now defined by:

$$\frac{\partial \gamma}{\partial x} = u\left(1+\delta\frac{\partial^2 H}{\partial p \partial x}\right) - v\,\delta\left(\frac{\partial^2 H}{\partial x^2} + p\frac{\partial^2 H}{\partial z \partial x}\right) - w(\,\delta\left(\frac{\partial H}{\partial x} - p\frac{\partial^2 H}{\partial p \partial x}\right) - \frac{1}{x} = 0,$$
$$\frac{\partial \gamma}{\partial p} = u\delta\frac{\partial^2 H}{\partial p^2} + v - v\,\delta\left(\frac{\partial^2 H}{\partial p \partial x} + p\frac{\partial^2 H}{\partial z \partial p} + \frac{\partial H}{\partial z}\right) + w\delta p\frac{\partial^2 H}{\partial p^2} - \frac{1}{p} = 0,$$
$$\frac{\partial \gamma}{\partial z} = u\delta\frac{\partial^2 H}{\partial p \partial z} - v\,\delta\left(\frac{\partial^2 H}{\partial z \partial x} + p\frac{\partial^2 H}{\partial z^2}\right) + w(1-\,\delta\left(\frac{\partial H}{\partial z} - p\frac{\partial^2 H}{\partial p \partial z}\right) - \frac{1}{z} = 0,$$

$$\partial^2 \widehat{H} = \begin{bmatrix} \frac{\partial^2 H}{\partial p \partial x} & -\frac{\partial^2 H}{\partial x^2} - p\frac{\partial^2 H}{\partial z \partial x} & -\frac{\partial H}{\partial x} + p\frac{\partial^2 H}{\partial p \partial x} \\ \frac{\partial^2 H}{\partial p^2} & -\frac{\partial^2 H}{\partial p \partial x} - p\frac{\partial^2 H}{\partial z \partial p} - \frac{\partial H}{\partial z} & p\frac{\partial^2 H}{\partial p^2} \\ \frac{\partial^2 H}{\partial p \partial z} & -\frac{\partial^2 H}{\partial z \partial x} - p\frac{\partial^2 H}{\partial z^2} & -\frac{\partial H}{\partial z} + p\frac{\partial^2 H}{\partial p \partial z} \end{bmatrix}$$

Let $U = (u, v, w)$ and $\frac{1}{P} = (\frac{1}{x}, \frac{1}{p}, \frac{1}{z})$, the critical point is defined by:

$$U + \delta\partial^2 \widehat{H} U - 1/p = 0,$$

and we apply the asymptotic probabilistic method. Various problems in physics concerning objects in movement obey more or less directly the Hamilton's equations $H(p, q, t)$. They must obey asymptotically differential partial equations of second order. When the movements are autonomous, the equations are the Laplacians.

## Conclusion

After this study, we can say, under good conditions, that an EDO is deterministic near the origin of the process, but may have stochastic or fixed cycles after a very long time.
With this probabilistic method, we obtain some new results, but we meet also many new difficulties due to the particular steepest descent's method used to study the Plancherel-Rotach's function.